\documentclass[11pt,a4paper,reqno]{amsart}
\usepackage{amssymb}
\usepackage{amscd}
\usepackage{color}
\usepackage{hyperref}
\numberwithin{equation}{section}

\numberwithin{equation}{section}

\theoremstyle{plain}

\newtheorem{theorem}[subsection]{Theorem}
\newtheorem{proposition}[subsection]{Proposition}
\newtheorem{lemma}[subsection]{Lemma}

\newtheorem{remark}[subsection]{Remark}

\theoremstyle{definition}

\newtheorem{definition}[subsection]{Definition}
\newtheorem{remarks}[subsection]{Remarks}
\newtheorem{example}[subsection]{Example}
\newtheorem{examples}[subsection]{Examples}

\renewcommand{\leq}{\leqslant}
\renewcommand{\geq}{\geqslant}

\newsavebox{\proofbox}
\savebox{\proofbox}{\begin{picture}(7,7)%
  \put(0,0){\framebox(7,7){}}\end{picture}}

\newcommand\Z{\mathbb{Z}}
\newcommand\R{\mathbb{R}}
\newcommand\T{\mathbb{T}}
\newcommand\C{\mathbb{C}}
\newcommand\N{\mathbb{N}}
\newcommand\D{\mathcal{D}}

\def\proof{\noindent\textit{Proof. }}

\def\remark{\noindent\textit{Remark. }}

\def\endproof{\hfill{\usebox{\proofbox}}}

\begin{document}

\title[Long-time instability and unbounded orbits]{Long-time instability and unbounded Sobolev orbits for some periodic nonlinear Schr\"odinger equations}
\author{Zaher Hani}

\address{Courant Institute of Mathematical Sciences, 251 Mercer Street, New York, NY 10012.}
\email{hani@cims.nyu.edu}
\thanks{The author is partly supported by a Simons Postdoctoral Fellowship.}

\maketitle

\begin{abstract}
We study the energy cascade problematic for some nonlinear Schr\"odinger equations on $\T^2$ in terms of the growth of Sobolev norms. We define the notion of long-time strong instability and establish its connection to the existence of unbounded Sobolev orbits. This connection is then explored for a family of cubic Schr\"odinger nonlinearities that are equal or closely related to the standard polynomial one $|u|^2u$. Most notably, we prove the existence of unbounded Sobolev orbits for a family of Hamiltonian cubic nonlinearities that includes the resonant cubic NLS equation (a.k.a. the first Birkhoff normal form).

\end{abstract}

\section{Introduction}\label{Introduction}

This manuscript is concerned with the dynamical behavior of solutions of the space-periodic nonlinear Schr\"odinger equation:
\begin{equation}\label{NLS}
\begin{cases}
-i\partial_t u +\Delta u&=\mathcal{N}(u)\\
u(0,x)&= u_0(x), \quad x\in \T^2
\end{cases}
\end{equation}
where $\mathcal N(u)$ stands for a cubic Hamiltonian nonlinearity (i.e. $\mathcal N(u)=\frac{\partial}{\partial \overline u}G(u,\overline u)$) that is either equal or closely related to the polynomial cubic one given by $\mathcal {N}_{\infty}(u)=|u|^2u$. We will work in spatial dimension $d=2$, but all the results extend directly to higher dimensions by considering solutions that depend only on the first two spatial coordinates. The initial data will be assumed to belong to some Sobolev space $H^s(\T^2)$ with $s>0$. 

We are particularly interested in the broad question of whether there exists global solutions to $\eqref{NLS}$ that exhibit oscillations at smaller and smaller spatial scales as time becomes larger and larger. In frequency space, this corresponds to solutions whose \emph{conserved} energy moves from low-frequency concentration zones to high-frequency ones, in what is commonly referred to as the \emph{forward or direct cascade}. One way to capture this energy dynamic in frequency space is to look at the behavior of Sobolev norms given by:
\begin{equation*}
\|u(t)\|_{H^s(\T^2)}=\left(\sum_{n\in \Z^2} \langle n \rangle^{2s}|\widehat{u}(t,n)|^2\right)^{1/2}
\end{equation*}
whose growth in time, for $s>1$, signals that a larger contribution to $\|u(t)\|_{H^s}$ is coming from higher and higher frequencies. The main question in this perspective is to find global solutions $u(t)$ of \eqref{NLS} for which:
\begin{equation}\label{unbdd orbit}
\sup_{t\in \R}\|u(t)\|_{H^s(\T^2)}=+\infty?
\end{equation}

A different perspective to energy cascade is provided by the physical theory of weak turbulence, which can be, loosely speaking, defined as the out-of-equilibrium statistics of random nonlinear waves \cite{DNPZ,ZLF}. This theory was invigorated in the 60s of the past century in plasma physics, where the subject adopted its name, and in the study of water waves. In these contexts, energy cascades were used to explain particle and energy transfer across tokomak plasma or how momentum is transferred from wind to the ocean. In contrast to the somewhat deterministic nature of question \eqref{unbdd orbit} above, weak turbulence theory addresses the cascade issue by looking at wave spectra (the average energy or mass density as a function of the wavelength) rather than individual wave trajectories. These spectra turn out to satisfy (after appropriate averaging and limiting operations) certain kinetic equations, for which Zakharov was able to find particular solutions (now called Kolmogorov-Zakharov spectra\footnote{Because they resemble Kolmogorov spectra in the theory of hydrodynamic turbulence.}) that correspond to constant energy flux through scales \cite{Zakh65}. Of course, a constant energy flux through scales is often not consistent with the conservation of energy of Hamiltonian systems like \eqref{NLS}, a fact which necessitates the presence of feeding source and a dissipative sink of energy located at large and small scales respectively.

While these spectra offer a valid evidence (from a physics point of view) of some form of cascade phenomenon happening at the medium (unexcited) frequencies, many of the manipulations involved are still far from rigorous. In this paper, we will adopt the first perspective mentioned above of searching for unbounded Sobolev orbits as in \eqref{unbdd orbit}. This seems to ask for a somewhat stronger form of cascade in comparison to the second perspective, which might explain why a satisfactory answer to \eqref{unbdd orbit} is still open for most equations of interest, including cubic NLS. For the polynomial cubic nonlinearity $\mathcal{N}_\infty(u)=|u|^2u$, the conservation laws of mass given by:

\begin{equation}\label{conservation of mass c7}
M[u(t)]:=\int_{\T^2} |u(t,x)|^2 dx=M[u(0)]
\end{equation}
and energy given by 
\begin{equation}\label{conservation of energy}
E[u(t)]:=\int_{\T^2}\frac{1}{2}|\nabla u(t,x)|^2+\frac{1}{4}|u(t,x)|^4 dx=E[u(0)]
\end{equation}
give global \emph{a priori} bounds on the $L^2$ and $H^1$ norms of the solution. So the answer to the question in \eqref{unbdd orbit} is negative for $s=0,1$. We should note that the one-dimensional analogue of the cubic NLS equation is completely integrable \cite{ZS} and higher conservation laws give global a priori bounds for $H^s(\T^2)$ at least when $s\in \N$. However, it is expected that this behavior is particular to the one-dimensional case and that, for $s>0, s\neq 1$, solutions exhibiting unbounded orbits in $H^s(\T^2)$ actually exist and are abundant. Proving this seems to be quite a difficult task that was considered by Bourgain in \cite{B6} as one of the next century problems in Hamiltonian PDE.

It is worth mentioning that polynomial-in-time upper bounds on the growth of $\|u(t)\|_{H^s}$ have been proved by several authors starting with pioneering works of Bourgain \cite{B7} and Staffilani \cite{Staffilani} (see also \cite{CDKS,CKSTT5, B10} and related work in \cite{CW, Zhong, Soh, CKO}). But this polynomial growth is expected to be far from optimal. In fact, Bourgain conjectured in \cite{B7} that the growth should be sub-polynomial in the sense that:
$$
\|u(t)\|_{H^s} \leq C_\epsilon t^\epsilon \quad \hbox{ for any } \epsilon>0,
$$
but this remains open as well.

We end this background overview with a remark concerning the effect of the growth of high Sobolev norms on the frequency dynamics of solutions. The conservation of mass and energy imply that any increase in high Sobolev norms is not only accompanied by a migration of energy to high frequencies, but also a migration of mass to low frequencies. This phenomenon is known in the physics literature as the \emph{inverse cascade of mass}, and it balances the forward cascade of energy to keep the system conservative\footnote{This is the exact analogue of the forward and backward cascade of enstrophy and energy (respectively) discovered by Kraichman \cite{Krai} in the theory of hydrodynamic turbulence.}.

\subsection{Long-time strong instability} Despite attracting a lot of attention, the existence of unbounded orbits as in \eqref{unbdd orbit} remains unproved for the polynomial cubic nonlinearity. However, some results exhibiting some form of cascade that is strictly weaker than $\eqref{unbdd orbit}$ have been obtained in the past years (see for example \cite{Kuksin, CKSTT2, CF, GuaKa} and \cite{GG1, GG2, Poc1, Poc2, Gua} for work on other Hamiltonian PDE). Of particular importance for this article is the result of Colliander, Keel, Staffilani, Takaoka, and Tao in \cite{CKSTT2} (see also \cite{GuaKa} for a refinement) who proved that for any given parameters\footnote{The theorem is stated in \cite{CKSTT2} for $s>1$. However, minor modifications to the proof allow one to extend that range to $s>0, s\neq 1$. See \cite{Thesis} for details.}
$s\in (0, \infty)\setminus \{1\}$, $K\gg 1$, and $0<\delta \ll 1$, there exists a global solution $u(t,x)$ of \eqref{cubic NLS} and a time $T>0$ such that
\begin{equation}
\label{CKSTT}
\tag{CKSTT}
\|u(0)\|_{H^s}\leq \delta \quad \hbox{and}\quad \|u(T)\|_{H^s}\geq K\, .
\end{equation}

The above result is equivalent to saying that the orbit of any $H^s$-\emph{neighborhood} of the origin under the nonlinear flow of the (polynomial) cubic NLS is not uniformly bounded in $H^s$. This is still considerably weaker than $\eqref{unbdd orbit}$ that requires a single orbit to be unbounded. 

It will be useful to interpret the result \eqref{CKSTT} above as the \emph{long-time strong instability of the flow near 0}, a concept that can be defined in great generality as follows:
\begin{definition}[Long-time strong instability] 
Let $X$ be a Banach space and suppose that $S(t)$ is a continuous dynamical system\footnote{This means that for every $x_0\in X$, $S(t)x_0 \in C_t([-T,T];X)$ for any $T>0$, and satisfies $S(0)x_0=x_0$ and $S(t_1+t_2)x_0=S(t_1)S(t_2)x_0$ for any $t_1, t_2 \in \R$.} on $X$. We say that the flow exhibits \emph{long-time strong instability near $\phi \in X$} if for every $0<\delta<1$ and $K\geq 1$, there exists $\phi^* \in B_X(\phi , \delta)$ such that $\sup_{t>0} \|S(t) \phi^*\|_X \geq K$.
\end{definition}

While the long-time strong instability near some $\phi$ is much weaker than the existence of an unbounded orbit in $X$ satisfying $\sup_{t}\|S(t) \phi\|_X=\infty$, it turns out that a sufficient extension of it is enough to guarantee the \emph{abundance} of unbounded orbits. We state this observation in the proposition below also in the generality of continuous dynamical systems on Banach spaces. For this, we will need an additional assumption on the flow, which is \emph{well-posedness}: We say that the flow $S(t)$ is well-posed if for every $T>0$, the map
$$
\phi \mapsto S(t)\phi$$
is continuous from $X \to \in C_ t([-T,T]\times X)$.

\begin{proposition}\label{observation} \textbf{\emph{(From long-time strong instability to generic unbounded orbits)}}
Let $X$ be a Banach space and $S(t)$ be a well-posed continuous dynamical system on $X$.
\begin{enumerate}

\item If the flow $S(t)$ exhibits long-time strong instability near a dense set of initial data $\phi_0 \in X$, then the set of unbounded orbits is generic\footnote{I.e. the intersection of countably many dense open sets.} in the sense of Baire category; in particular it is nonempty. More precisely, the set
$$
\{\phi_0 \in X: \sup_{t>0} \|S(t) \phi_0\| <\infty\} \quad \hbox{is meager (countable union of nowhere dense sets)}. 
$$
\medskip
\item Let $\mathcal{F}$ be a closed subset of $X$ and $D \subset \mathcal F$ be dense in $\mathcal F$. Suppose that for any $\phi \in D, \delta>0, K\geq 1$, there exists $\phi^* \in \mathcal F$ such that $\|\phi-\phi^*\|_{X}\leq \delta$ and $\sup_{t>0}\|S(t)\phi^*\|_{X} \geq K$. Then the set $\{\psi \in \mathcal F: \sup_{t\in \R} \|S(t)\psi\|_X=\infty\}$ is co-meager in $\mathcal{F}$ (in particular non-empty).

\end{enumerate}
\end{proposition}

The above proposition can somehow be understood as a ``nonlinear uniform boundedness principle" and its proof is a straightforward application of the Baire Category theorem. Indeed, let

$$\mathcal{B}=\{u_0 \in \mathcal{F}: \sup_{t>0}\|S(t)u_0\|_{X}<\infty\}$$
and 
$$\mathcal{B}_N=\{u_0 \in \mathcal{F}: \sup_{t>0}\|S(t)u_0\|_{X}\leq N, \}.$$
Obviously, $\mathcal{B}=\cup_{N\in \Z, N\geq 1}\mathcal{B}_N$. $\mathcal{B}_N$ is closed in $X$ courtesy of the well-posedness of the flow $S(t)$. The long-time strong instability of the flow near every $u_0 \in D$ and the fact that $D$ is dense in $\mathcal F$ implies that $\mathcal{B}_N$ has no interior as a subset of $\mathcal{F}$ (because for any $u_0\in D $ and any $\delta>0$, one can find $v_0 \in B_X(v_0, \delta)\cap \mathcal{F}$ whose orbit exits the ball $B_{X}(0,2N)$). As a result, $\mathcal{B}\subset \mathcal{F} $ is a countable union of closed no-where dense sets. Since $\mathcal{F}$ is a complete metric space, $\mathcal{F}\setminus\mathcal{B} \neq \emptyset$ and the set of unbounded orbits is co-meager in $\mathcal{F}$.

In practice, proving the genericness of unbounded orbits by establishing long-time strong instability near a dense set of $X$, as in part (i) of the above proposition, can be quite a demanding task. Nonetheless, the existence of \emph{at least one} unbounded orbit can be seen to trivially imply the validity of part (ii) of the proposition for \emph{some} closed subset of $X$, which makes the implication statement in that part an equivalence. In what follows, we will apply the above proposition for the NLS flow in \eqref{NLS} and prove long-time strong instability near some carefully chosen families of $\phi \in X$. This will be sufficient to invoke part (ii) of Proposition \ref{observation} in order to conclude the existence of unbounded orbits for certain cubic NLS nonlinearities that are arbitrarily close (but not quite equal) to the polynomial one.

\subsection{Implications for NLS} Proposition \ref{observation} suggests a program to proving the existence of unbounded orbits of \eqref{NLS} by proving analogues of the \eqref{CKSTT} instability result in \cite{CKSTT2} near sufficiently large families of initial data in $H^s$. Indeed, with the Banach space $X$ taken to be $H^s(\T^2)$, the NLS flow defines a well-posed continuous dynamical system whenever equation \eqref{NLS} is globally well-posed in $H^s(\T^2)$\footnote{For instance, this is the case for the polynomial cubic nonlinearity whenever $s>2/3$ \cite{B4, dSPST}.}. 

For the polynomial cubic nonlinearity $|u|^2u$, we take a rather timid step in the direction suggested by Proposition \ref{observation}, and prove long-time strong instability of the flow in $H^s(\T^2)$ near plane wave initial data given by $u_0^*=Ae^{inx}$ ($A\in \C$, $n\in \Z^2$) in the range $0<s<1$:

\begin{theorem}\label{plane wave instability}
Consider the polynomial cubic NLS equation:
\begin{equation}\label{cubic NLS}
\begin{cases}
-i\partial_t u +\Delta u&=|u|^2u\\
u(0,x)&= u_0(x)\in H^s(\T^2)
\end{cases}
\end{equation}
For any $0<s<1$, $0<\delta \ll1$, $K\gg1$, there exists a global solution $u(t)$ such that $\|u(0)-Ae^{inx}\|_{H^s(\T^2)}\leq \delta$ and $\|u(T)\|_{H^s(\T^2)} \geq K$ for some later time $T$.
\end{theorem}

Solutions to \eqref{cubic NLS} with initial data $Ae^{in.x}$ are explicit and are given by $Ae^{i(n.x +|n|^2t+|A|^2t)}$. While the proof of this theorem relies heavily on the constructions in \cite{CKSTT2}, it offers a stark example of some of the difficulties encountered as one tries to prove long-time strong instability for \eqref{cubic NLS} near non-zero initial data. In fact, the major issue in proving Theorem \ref{plane wave instability} is not only to find a tractable mechanism for norm growth (which will be a modification of that in \cite{CKSTT2}), but also to find ways to limit or incorporate interactions between two different parts of the solution: one is the ``base solution" that has $O(1)$ size in $H^s$ and is a manifestation of the unperturbed plane wave solution $Ae^{i(n.x +|n|^2t+|A|^2t)}$; and the other is the ``turbulent" part that has $O(\delta)$ size in $L^2$, but whose Sobolev norm has to grow to become larger than $K$ after time $T$. 

While we believe that the restriction in Theorem \ref{plane wave instability} to the range $0<s<1$ is probably technical, we should mention that for large enough $s>1$, plane wave solutions were shown in \cite{FGL} to exhibit polynomial-time Nekhorashev-type stability in the following sense: For any $N \in \N$, there exists $\epsilon_0$ such if $\epsilon<\epsilon_0$, most solutions that are initially $\epsilon$-close to a plane wave $Ae^{in.x}$ in $H^s$ (for $s$ large enough) will remain as such\footnote{I.e. they remain in a 2$\epsilon$-neighborhood of a plane wave $A(t)e^{in.x}$.} for a period of time that is $O(\epsilon^{-N})$. Of course, this does not prohibit growth of the Sobolev norm in infinite time, but it suggests the need for a much more elaborate growth mechanism that beats this long stability time interval. Indeed, it might be worth mentioning that the same interactions that require us to restrict the instability Theorem \ref{plane wave instability} to the range $0<s<1$ are the ones that allow the authors in \cite{FGL} to prove their Nekhorashev-type stability result for large enough $s$.

While we still cannot push the program suggested by Proposition \ref{observation} to its end for the polynomial cubic nonlinearity, a positive answer to question \eqref{unbdd orbit} can be given for some cubic nonlinearities that are ``almost polynomial". By this we mean that they are obtained from the polynomial cubic nonlinearity $|u|^2u$ by eliminating \emph{highly non-resonant}\footnote{Non-resonant interactions have a far weaker effect on the dynamical behavior than resonant ones (at least for early time scales), which makes the result in Theorem \ref{WT for special nonlinearities} further evidence of the existence of unbounded orbits for the polynomial cubic equation \eqref{cubic NLS}.} terms beyond any arbitrary threshold.

More precisely, our second result proves the existence of unbounded orbits for the following family of \emph{Hamiltonian} nonlinearities:
\begin{equation}\label{def of N_R}
\mathcal{N}_R(u)=\sum_{n \in \Z^2} \left(\sum_{\Gamma_R(n)} \widehat{u}(n_1) \overline{\widehat{u}(n_2)} \widehat u(n_3)\right) e^{2\pi i n.x};\quad R\in \{0\}\cup \N.
\end{equation}
where $\Gamma_R(n)=\{ (n_1, n_2, n_3) \in (\Z^2)^3: n_1-n_2+n_3=n \hbox{ and } \omega_4:= |n_1|^2-|n_2|^2+|n_3|^2-|n|^2 \in [-R, R]\}$.

Note that the polynomial cubic nonlinearity $\mathcal{N}_\infty(u)=|u|^2u$ corresponds to taking $R=\infty$ in the above definition. Another well-known nonlinearity that is included in this family is the resonant cubic NLS system (a.k.a. the first Birkhoff normal form) that corresponds to $R=0$. We remark that for any $R\in \{0\}\cup \N \cup \{\infty\}$, the nonlinearity $\mathcal{N}_R(u)$ is Hamiltonian and conserves mass. In fact, in the two extreme cases $R=0$ and $R=\infty$, the Hamiltonian is positive definite and controls the $\dot H^1$ norm of the solution (cf. \eqref{conservation of energy} and \eqref{R=0 Hamiltonian}). The same is true for any $R\in \N$ provided that the mass of the solution is smaller than some positive threshold (cf. Section \ref{LT instability section}). This will allow us to conclude global well-posedness of the flow (cf. Section \ref{LT instability section}) in a way that allows to apply the paradigm of Proposition \ref{observation} and prove our main result:

\begin{theorem}\label{WT for special nonlinearities}
For any $R\in \{0\}\cup \N$ and any $s>1$, there exists global solutions to \eqref{NLS} (with $\mathcal{N}(u)=\mathcal{N}_R(u)$) whose orbit is unbounded in $H^s(\T^2)$ as in \eqref{unbdd orbit}. Moreover, initial data corresponding to unbounded orbits are residual in a closed subset of $H^s(\T^2)$.
\end{theorem}

A few remarks are in order: First, the same conclusion holds in the range of regularities $0<s<1$ under the conditional assumption that the system is globally well-posed for sufficiently small data (cf. Section \ref{LT instability section}). Second, the closed set $\mathcal F \subset H^s(\T^2)$ mentioned in Theorem \ref{WT for special nonlinearities}, that contains the residual set of unbounded orbits, is defined in terms of a geometric and combinatorial condition imposed on the Fourier support. This condition allows one to prove long-time strong instability near a dense subset and conclude as in part (ii) of Proposition \ref{observation} that bounded orbits are a meager subset of $\mathcal F$. Finally, we should mention that Bourgain constructed in \cite{B11} a nonlinearity $\mathcal N(u)$ that grows linearly in $u$ for which he proved the existence of an unbounded orbit. The novelty in Theorem \ref{WT for special nonlinearities} comes from 1) proving the \emph{abundance} of unbounded orbits for a family of Hamiltonian \emph{cubic} nonlinearities $\mathcal N_R(u)$ that approximates the polynomial cubic one ($R=\infty$) and includes the well-known resonant cubic system ($R=0$); and 2) illustrating how the program suggested by Proposition \ref{observation} can be applied to conclude existence and genericness (possibly in some metric subspace) of unbounded orbits.

The paper is organized as follows: in Section \ref{Resonant and finite reductions}, we set up our problem in Fourier space and present some basic reductions and notions that will be used in the rest of the paper. In Section \ref{Proof of instability section}, we prove Theorem \ref{plane wave instability}; and finally Theorem \ref{WT for special nonlinearities} is proved in Section \ref{LT instability section}. Throughout this paper, We write $A\lesssim B$ to signify that there is a constant $C>0$ such that $A\le C B$. We also write $A\sim B$ when $A\lesssim B\lesssim A$. If the constant $C$ involved has some explicit dependency, we emphasize it by a subscript. Thus $A\lesssim_uB$ means that $A\le C(u)B$ for some constant $C(u)$ depending on $u$. In some instances, we use the notation $A \ll B$ to signify that the implicit constant $C$ is large. Finally, we omit all factors of $2\pi$ arising in the definition of the Fourier transform on $\T^2$ as they play no role in our analysis.

\subsection*{Acknowledgements} The author would like to thank Terry Tao for introducing him to this line of research as a Ph.D. student, and for all the invaluable discussions and suggestions. The author is supported by a Simons Postdoctoral Fellowship.

\section{Resonant and finite reductions}\label{Resonant and finite reductions}
We start by considering the equation:
\begin{equation}\label{cubic NLS}
\begin{cases}
-i\partial_t u +\Delta u&=|u|^2u\\
u(0,x)&= u_0(x)
\end{cases}
\end{equation}

In this section, we perform some reductions that will be used in the course of the paper. First, we will cast \eqref{cubic NLS} as an infinite system of ODE which we call $(\mathcal{F}NLS)$ in $\{a_n(t)\}_{n\in \Z^2}$ where $a_n(t)$ are closely related to the Fourier coefficients of the solution $u(t,x)$. Afterwards we isolate the resonant reduction or the first Birkhoff normal form of this system $(\mathcal{RF}NLS)$, which is obtained from $(\mathcal{F}NLS)$ by deleting all non-resonant interactions. We then explain how to reduce further this resonant system into a finite ``Toy system" of ODE as was done in \cite{CKSTT2} by choosing initial data supported on special subsets of $\Z^2$. These reductions, borrowed from \cite{CKSTT2} and summarized here for completeness, form an important ingredient for the proof of Theorems \ref{plane wave instability} and \ref{WT for special nonlinearities}.

\subsection{NLS as an infinite system of ODE $\mathcal{F}NLS$}\label{reductions subsection}

Equation $\eqref{cubic NLS}$ enjoys the following ``guage freedom": 
$$
v(t,x)=e^{iGt}u(t,x)
$$ 
where $v$ satisfies the equation:
$$
-i\partial_t v+\Delta v=(G+|v|^2)v
$$
Using the ansatz:
$$
v(t,x)=\sum_{n}a_n(t)e^{i(n.x+|n|^2t)}
$$
we get after a direct calculation that the coefficients $a_n(t)$ satisfy the following system of ODE:

\begin{equation}\label{a_n ODE2}
-i\partial_t a_n(t)=(G+2M)a_n-|a_n|^2a_n+\sum_{\substack{n_1-n_2+n_3=n\\n_1,n_3\neq n}}a_{n_1}\overline{a_{n_2}}a_{n_3}e^{i\omega t}
\end{equation}
where $M$ stands for the mass :
$$
M=\sum_n |a_n|^2
$$
which is conserved in time. Choosing $G=-2M$ we end up with an infinite system of ODE which we denote by $(\mathcal{FNLS})$ given by:

\begin{equation}\label{FNLS}
-i\partial_t a_n=-a_n|a_n|^2+\sum_{n_1,n_2,n_3 \in \Gamma(n)}a_{n_1}\overline{a_{n_2}}a_{n_3}e^{i\omega_4 t}
\end{equation}
where $\Gamma(n)=\{(n_1,n_2,n_3)\in (\Z^2)^3: n_1-n_2+n_3=n \text{ and } n_1,n_3 \neq n\}$.

We remark that at any time $t$, one can recover $u$ from $v$ and vice-versa and the two solutions have the same Sobolev norms. It is easy to see that this system of ODE is locally well-posed in $l^1(\Z^2)$ (see \cite{CKSTT2}).

\subsection{Resonant reduction $\mathcal{RF}NLS$} 
The first interaction term in $\eqref{FNLS}$ is a self interaction term, whereas the second term is a sum taken over all parallelograms (which may collapse to segments) with $n$ as one of its vertices (and $n_2$ being the opposite vertex). Analyzing the resonance factor $\omega_4$ we get that the parallelogram is a rectangle exactly when $\omega_4=0$. In fact,

\begin{align*}
w_4=&|n_1|^2-|n_2|^2+|n_3|^2-|n|^2= |n_1-n|^2-|n_2-n|^2+|n_3-n|^2\\=&-2(n_1-n).(n_3-n).
\end{align*}

The interactions for which $\omega_4=0$ are called \emph{resonant interaction}, whereas those for which $\omega_4\neq 0$ are \emph{non-resonant}. Resonant interactions will play a particularly important role in the analysis, a fact which motivates us to define the set of all \emph{resonant non-self interactions} $\Gamma_{res}(n)\subset \Gamma(n)$ as:

\begin{equation}\label{def of Gamma_res}
\Gamma_{res}(n)=\{(n_1,n_2,n_3)\in \Gamma(n): \omega_4=|n_1|^2-|n_2|^2+|n_3|^2-|n|^2=0\}.
\end{equation}

The resonant reduction $(\mathcal{RF}NLS)$ of $(\mathcal{F}NLS)$ is obtained by deleting all non-resonant interactions and retaining the resonant ones only. More precisely, $(\mathcal{RF}NLS)$ is the system:

\begin{equation}\label{RFNLS}
-i\partial_t r_n=-r_n|r_n|^2+\sum_{n_1,n_2,n_3 \in \Gamma_{res}(n)}r_{n_1}\overline{r_{n_2}}r_{n_3}
\end{equation}

This system can be also obtained as the first Birkhoff normal form of \eqref{FNLS}. It will be important to notice that all non-self interactions (the second term on the RHS of $\eqref{RFNLS}$) of the resonant system $(\mathcal{RF}NLS)$ occur in rectangles: each $n\in \Z^2$ interacts with other frequencies by forming a rectangle with $n$ as one of its vertices. This will be the key to collapsing this system to a finite system of ODE by choosing the initial data to be supported on special subsets of $\Z^2$ that are somewhat \emph{``closed under completing rectangles"}. This is the content of the first two properties below.

\subsection{Nuclear families and the finite reduction \cite{CKSTT2}}\label{Nuclear families}

Suppose that for some subset $\Lambda \subset \Z^2$ we have:

\begin{itemize}
\item Property I$_\Lambda$ (Initial data): The initial data $r_n(0)$ is supported on $\Lambda$, i.e. $r_n(0)=0$ unless $n\in \Lambda$.
\item Property II$_\Lambda$ (Closure): If $n_1,n_2, n_3 \in \Lambda$ are three vertices of a rectangle, then the forth vertex of that rectangle is also in $\Lambda$. This implies that if $(n_1,n_2,n_3) \in \Gamma_{res}(n)\cap \Lambda^3$, then $n\in \Lambda$.
\end{itemize}

These two properties imply that the support of $r_n(t)$ remains restricted to $\Lambda$ for all later times. In fact, if $B(t):=\sum_{n\notin \Lambda}|r_n(t)|^2$, then it is easy to see that $|B'(t)|\leq C|B(t)|$ and the claim follows by applying Gronwall's inequality (see \cite{CKSTT2} for details or Lemma \ref{Gronwall on S} for a similar argument).

This observation allows to reduce the infinite system of ODE $(\mathcal{RF}NLS)$ into a finite system of ODE by restricting the initial data to a finite subset $\Lambda$ that satisfies the closure property above. As was done in \cite{CKSTT2}, we will further reduce this finite system into a smaller one by imposing some extra structure on the set $\Lambda$. This structure comes from the following hierarchy:

We will assume that $\Lambda$ splits into a disjoint union $\Lambda=\Lambda_1 \cup \ldots \cup \Lambda_P$ for some integer $P>1$ and call $\Lambda_j$ the \emph{$j-$th generation} of $\Lambda$. We define a \emph{nuclear family} in $\Lambda$ as a rectangle $(n_1, n_2, n_3,n_4)$ where the frequencies $n_1,n_3 $ (known as the ``parents") belong to the $j-$th generation $\Lambda_j$ while the ``\emph{children}" $n_2, n_4$ belong to $\Lambda_{j+1}$. Of course, if $(n_1,n_2,n_3,n_4)$ is a nuclear family then so are the \emph{trivial permutations} $(n_3,n_2,n_1,n_4)$ ,$(n_1,n_4,n_3,n_2)$, and $(n_3,n_4,n_1,n_2)$. In addition to the properties I$_\Lambda$ and II$_{\Lambda}$ above, we require:

\begin{itemize}
\item Property III$_{\Lambda}$ (Existence and uniqueness of spouse and children): For each $1\leq j <P$ and every $n_1\in \Lambda_j$, there exists a unique spouse $n_3\in \Lambda_j$ and unique (up to trivial permutations) children $n_2,n_4 \in \Lambda_{j+1}$ such that $(n_1,n_2,n_3,n_4)$ is a nuclear family in $\Lambda$.

\item Property IV$_{\Lambda}$ (Existence and uniqueness of parents and siblings): For each $1\leq j <P$ and evey $n_2 \in \Lambda_{j+1}$ there exists a unique sibling $n_4\in \Lambda_{j+1}$ and unique (up to permutation) parents $n_1,n_3 \in \Lambda_j$ such that $(n_1,n_2,n_3,n_4)$ is a nuclear family in $\Lambda$.

\item Property V$_\Lambda$ (Non-degeneracy): A sibling of any frequency $n$ is never equal to its spouse.
\item Property VI$_\Lambda$ (Faithfulness): Apart from nuclear families $\Lambda$ contains no other rectangles. In fact, by the closure property II$_\Lambda$, this also means that it contains no right angled triangles other than those coming from vertices of nuclear families.
\end{itemize}
Assuming that such a set exists, the system $(\mathcal{RF}NLS)$ can be written as:

$$
-i\partial_t r_n=-r_n|r_n|^2 +2r_{n_{child-1}}r_{n_{child-2}}\overline{r_{n_{spouse}}}+2r_{n_{parent-1}}r_{n_{parent-2}}\overline{r_{n_{sibling}}}
$$
where we used the fact that the only rectangles that contain $n$ as a vertex are those coming from nuclear families where $n$ is either a parent (second term on the RHS) or a child (third term on RHS). The constant factor 2 in front of these terms comes from the trivial permutations of the nuclear families. 

The final property imposed on the initial data allows to collapse the $(\mathcal{RF}NLS)$ system further.
\begin{itemize}
\item Property VII$_\Lambda$ (Intragenerational equality): The function $n\mapsto r_n(0)$ is constant on each generation $\Lambda_j$. In other words, $r_{n_1}(0)=r_{n_2}(0)$ whenever $n_1,n_2 \in \Lambda_j$.
\end{itemize}

One can verify (either by a Gronwall or a bootstrap argument) that Property VII$_\Lambda$ is preserved by the flow of $(\mathcal{RF}NLS)$ (cf. \cite{CKSTT2}) in the sense that intragenerational equality remains true for all time. By setting $r_n(t)=b_j(t)$ whenever $n \in \Lambda_j$, we arrive at the following system for $b_j$:

\begin{equation}\label{toy system}
-i\partial_t b_j=-b_j|b_j|^2+2b_{j+1}^2\overline{b_j}+2b_{j-1}^2\overline{b_j}
\end{equation}
with the convention\footnote{Actually, one can easily verify that if one regards $\eqref{toy system}$ as a system of ODE in $b_j$ with $j\in \Z$, then the property $b_j(0)=0$ for $j\notin \{1,\ldots, P\}$ is propagated by the flow.} that $b_0(t)=b_{P+1}(t)=0$.

With these reductions at hand, the result \eqref{CKSTT} was proved in \cite{CKSTT2} in three main steps:

\begin{enumerate}
\item Construction of the frequency set $\Lambda$ satisfying Properties I$_{\Lambda}$-VII$_{\Lambda}$ above, along with the norm explosion property:

\begin{equation}\label{norm explosion}
\frac{\sum_{n \in \Lambda_{P-2}}|n|^{2s}}{\sum_{n\in \Lambda_3}|n|^{2s}} \gtrsim \frac{K^2}{\delta^2}
\end{equation}

\item Construction of a solution of resonant system \eqref{RFNLS} that is initially concentrated on a low frequency generation (namely $\Lambda_3$ for technical reasons) and after some time $T$ concentrated on a high-frequency generation $\Lambda_{N-2}$. This amounts to proving the following property of the Toy system \eqref{toy system}:

\begin{proposition}\label{diffusion result}
Given $P \gg 1, \epsilon \ll 1$, there exists initial data $b(0)=(b_1(0),$ $\ldots, b_N(0))$ for $\eqref{toy system}$ and a time $T=T(P,\epsilon)$ so that 

\begin{align}
|b_3(0)|\geq 1-\epsilon,& \quad \hbox{ while }|b_j(0)|\leq \epsilon \hbox{ for } j\neq 3 \label{toy system t=0};\\
|b_{P-2}(T)|\geq 1-\epsilon,&\quad \hbox{ while } |b_j(T)|\leq \epsilon\hbox{ for } j\neq P-2. \label{toy system t=T}
\end{align}

Also, the solution satisfies $\|b(t)\|_{l^\infty}\sim 1$ for all $0\leq t\leq T$.
\end{proposition}
\end{enumerate}
The above two steps allow to prove the needed Sobolev norm growth for the resonant system $(\mathcal{RF}NLS)$ in \eqref{RFNLS}. To transfer this growth to $(\mathcal{F}NLS)$, one needs:
\begin{enumerate}
\setcounter{enumi}{2}
\item \textbf{Stability lemma:} In this step, one proves the faithfulness of the approximation of \eqref{FNLS} by \eqref{RFNLS} over a time interval that is \emph{strictly longer} than that needed for the Sobolev norm of \eqref{RFNLS} to grow. 
\end{enumerate}

\section{Long-time instability of plane waves}\label{Proof of instability section}

In this section, we prove Theorem \ref{plane wave instability}. We start by reducing to the case when $v=0$ by resorting to the Galilean symmetry: 
\begin{equation}\label{Galilean invariance}
u(t,x)\leftrightarrow e^{i(v.x+|v|^2t)}u(t,x+2vt)
\end{equation}
Indeed, if $\tilde u$ solves \eqref{cubic NLS}, $\|\tilde u(0)-A\|_{H^s(\T^2)}\leq \tilde \delta$, and $\|\tilde u(T)\|_{H^s(\T^2)} \geq \tilde K$, then $u(t)=e^{i(v.x+|v|^2t)}\tilde u(t,x+2vt)$ also solves \eqref{cubic NLS} and satisfies $\|u(0)-Ae^{iv.x}\|_{H^s(\T^2)}=\|e^{iv.x}(\tilde u(0)-A)\|_{H^s(\T^2)}\lesssim_{v} \|\tilde u(0)-A\|_{H^s(\T^2)} \lesssim \tilde \delta$, while $\|u(T)\|_{H^s(\T^2)}=\|e^{iv.x}\tilde u(T)\|_{H^s(\T^2)}\gtrsim_v \tilde K$. The result then follows by choosing $\tilde K$ and $\tilde \delta$ appropriately in terms of $K, \delta,$ and $v$. Also, we may assume without any loss of generality (thanks to phase translation invariance) that $A\in \R$ and positive.

The initial data will be taken to be $u(0)=A+w(0)$ where $w(0)$ is supported on a frequency set $\Lambda$ similar to that in the previous section. As before, the norm explosion will happen as the energy cascades between the different generations of $\Lambda$. For this to work, we need to limit the interactions between the large mode at $0$ and the modes in $\Lambda$, so that this cascade process is not disturbed or reversed. Some of these interactions can be eliminated by imposing an additional ``geometric" condition on $\Lambda$ (Property VIII$_{\Lambda}$ below). Others cannot be easily eliminated and they are the reason why we restrict our result to the regime $0<s<1$. More precisely, we are referring to non-resonant interactions that happen when a mode at $n \in \Lambda$ interacts non-resonantly with the zero mode to excite the mode at $-n$ which in turn interacts with the zero mode to excite the mode at $n$. The net self-excitation of the mode at $n$ due to this interaction chain is actually resonant (but at a second order) and leads to the restriction on $s$.

In contrast to the argument in \cite{CKSTT2}, it will be important to work on the dilated lattice $\mathcal D=N\Z^2$ generated by $(N,0)$ and $(0,N)$. It is easy to verify (via a Gronwall argument for example) that if the initial data is supported in frequency space on $\mathcal{D}$, then the emanating solution will also be supported on $\mathcal{D}$ for all time. This ultimately follows because the interactions of the $(\mathcal{F}NLS)$ system happen in parallelograms and $\mathcal{D}$ trivially satisfies the \emph{parallelogram closure property}: if three vertices of a parallelogram are in $\mathcal{D}$, then so is the forth.

Consequently, we will consider the $(\mathcal{F}NLS)$ system posed on $\mathcal D$:
\begin{equation}
-i\partial_t c_n(t)=-|c_n|^2c_n +\sum_{\Gamma(n)}c_{n_1}\overline{c_{n_2}}c_{n_3}e^{i\omega_4t}
\end{equation}
where $\Gamma(n)=\{(n_1,n_2,n_3)\in \D^3: n_1-n_2+n_3=n \text{ and }
n_1,n_3 \neq n\}$ and $\omega_4=|n_1|^2-|n_2|^2+|n_3|^2-|n|^2$. Notice that if $\omega_4 \neq 0$, then we directly get that $|\omega_4|\geq N^2$.  

Since we are considering solutions with a relatively large
zero-frequency mode, we have to isolate the terms of this sum
involving this mode. First, we isolate the $c_0(t)$ equation:

\begin{equation}\label{c_0}
-i\partial_t c_0(t)=-|c_0|^2c_0 +\overline{c_0(t)}\sum_{n_1\neq 0}c_{n_1}(t)c_{-n_1}(t)e^{2i|n_1|^2t}+\sum_{\tilde{\Gamma}(0)}c_{n_1}\overline{c_{n_2}}c_{n_3}e^{i\omega_4t}
\end{equation}
where $\tilde{\Gamma}(0)=\{(n_1,n_2,n_3)\in \D^3: n_1-n_2+n_3=0
\text{ and } n_1,n_2,n_3 \neq 0\}$

Another calculation shows that the equation of $c_n$ for $n\neq 0$ will have the following more elaborate form.

\begin{equation}\label{c_n}
\begin{split}
-i\partial_t c_n(t)=&-|c_n|^2c_n +c_0(t)^2\overline{c_{-n}(t)}e^{-2i|n|^2t}+2\left(\sum_{n_2\neq 0,-n}\overline{c_{n_2}}c_{n+n_2}e^{2i n\cdot n_2t}\right)c_0(t)\\
&+\left(\sum_{n_1\neq 0,n}c_{n_1}(t)c_{n-n_1}(t)e^{2in_1\cdot (n_1-n)t}\right)\overline{c_0}(t)+\sum_{\tilde{\Gamma}(n)}c_{n_1}\overline{c_{n_2}}c_{n_3}e^{i\omega_4t}
\end{split}
\end{equation}
where $\tilde{\Gamma}(n)=\{(n_1,n_2,n_3)\in \Gamma(n):,
n_1,n_2,n_3\neq 0\}$.

To construct a norm-exploding solution to $\eqref{c_n}$, we will first construct an energy-cascading solution of the resonant system $(\mathcal{RF}NLS)$ in a similar way to what was done in \cite{CKSTT2}. This solution to the resonant system will be blind to the large mode at $0$ and will be supported on a finite set $\Lambda$ located at high frequencies. The main difficulty is to show that on the set $\Lambda$, solutions to \eqref{c_n} can still be approximated by solutions to the resonant system in spite of the additional interaction terms with the large mode at 0.

Denote by $(\tilde c_n(t))$ the solution to the resonant system:
\begin{equation}\label{RFNLS v2}
-i\partial_t \tilde{c}_n(t)=-|\tilde{c}_n(t)|^2\tilde{c}_n(t)+\sum_{\Gamma_{res}(n)}\tilde{c}_{n_1}(t)\overline{\tilde{c}_{n_2}(t)}\tilde{c}_{n_3}(t)
\end{equation}
with suitably chosen initial conditions. $(\tilde c_n(0))$ will be finitely supported on a resonant set $\Lambda \subset \mathcal D\setminus 0$ that is closed under completing rectangles in the sense that it satisfies Property II$_\Lambda$ from the previous section. Recall that this implies that $(\tilde c_n(t))$ remains supported on $\Lambda$ for all later times. We will also impose the same conditions as before on the set $\Lambda$, namely Properties II$_\Lambda$-VI$_\Lambda$ as well as the norm explosion condition $\eqref{norm explosion}$. 

To be able to prove the closeness of $(\tilde c_n)$ and $(c_n)$ for $n \in \Lambda$, we will impose a new condition that is used to limit the interaction of the frequencies in $\Lambda$ with the 0-frequency mode $c_0$. 
\begin{itemize}

\item Property VIII$_\Lambda$ (No pair of modes from $\Lambda$ forms a rectangle with 0): $\Lambda$ contains at most one vertex of any rectangle having a vertex at 0, i.e. for any $n_1, n_2 \in \Lambda$, neither one of the two parallelograms determined by the three vertices $0,n_1,n_2$ is a rectangle. 

\end{itemize}

This property says that no resonant interactions can happen between two frequencies in $\Lambda$ and the zero mode. Such interactions with the relatively large zero mode can leak energy out of $\Lambda$ or inject energy into it, in a way that could destroy the delicate Arnold diffusion responsible for the energy cascade across the different generations of $\Lambda$.

The existence of such a $\Lambda$ is guaranteed by the following proposition:

\begin{proposition}\label{construction of Lambda v1}
Let $0<s, s\neq 1$ and $N\in \N$. Given parameters $\delta \ll 1$ and $K\gg1$, there exists $P=P(\delta, K)$  and a frequency set $\Lambda=\Lambda_1 \cup \ldots \cup \Lambda_P \subset \D=N\Z^2$ satisfying Property II$_\Lambda$-Property VIII$_\Lambda$ as well as the norm explosion condition $\eqref{norm explosion}$. In particular, the cardinality of $P$ can be chosen uniformly in $N$.
\end{proposition}

\proof: This will be a special case of Proposition \ref{construction of Lambda} in the next section (and the remark following its proof).
\endproof

The growth of the Sobolev norm of $(\tilde c_n)$ will follow by choosing the initial data as we recalled in the previous section and using the diffusion result in Proposition \ref{diffusion result}. The key point is then to transfer this growth property of $(\tilde c_n)$ into a statement about the Sobolev norm growth of $(c_n)$. This will be a consequence of the following lemma, which constitutes the bulk of the proof:

\begin{lemma}\label{stability lemma}
Fix $0<\sigma <1$ and suppose that $1\ll B$, $B^{1+\sigma}\ll N$, and $T\ll B^2\log B$ (implicit constants are allowed to depend on $\sigma$). Suppose that $\tilde c(t)=(\tilde c_n(t))$ is a solution to the resonant $\eqref{RFNLS v2}$ satisfying:
\begin{itemize}
\item $c_n(0)=0$ unless $n \in \Lambda$ for some finite set $\Lambda$ satisfying properties I$_\Lambda$-VIII$_\Lambda$ above.
\item \begin{equation}\label{bound on tilde c}
\|\tilde c\|_{l^1(\Z^2)}\leq B^{-1}\quad \hbox{ for all $0\leq t \leq T$.}
\end{equation}

\end{itemize}

If $c(t)$ is the solution to $(\mathcal{F}NLS)$ with initial data $c(0)=A+\tilde c(0)$ (i.e. $c_0(0)=A$ and $c_n(0)=\tilde c_n(0)$ for $n\neq 0$), then

\begin{equation}\label{stability estimate}
\| c(t)-\tilde c(t)\|_{l^1(\Lambda)}\lesssim_A \frac{B^{1+\sigma}}{N^2}
\end{equation}
for all $0 \leq t \leq T$.
\end{lemma}

\subsection{Proof of Theorem \ref{plane wave instability} assuming Lemma \ref{stability lemma}}\label{scaling argument}

Suppose that $0<s<1, \delta\ll1,$ and $K\gg 1$ are given. Choose $0<\sigma<1$ such that $s(1+\sigma)<1$ and let $P=P(\delta, K)$ and $\Lambda=\Lambda_1\cup \ldots \cup \Lambda_P$ be given by Proposition \ref{construction of Lambda v1}. Also let $b(t)=(b_1(t),\ldots, b_P(t))$ be the solution to the ``Toy system" $\eqref{toy system}$ over a time interval $[0,T_0]$ given by Proposition \ref{diffusion result} for some $\epsilon=\epsilon(\delta, K)$ to be specified later. Note that the ``Toy system" $\eqref{toy system}$ satisfies the following scaling symmetry:
$$
b^\lambda(t)=\frac{1}{\lambda}b(\frac{t}{\lambda^2})
$$
also solves $\eqref{toy system}$ on the interval $[0,\lambda^2T_0]$. Let $\tilde c(t)=(\tilde c_n(t))$ be the solution of the resonant system obtained from $b^\lambda(t)$ (i.e. $\tilde c_n(t)=b^\lambda_j(t)$ if $n\in \Lambda_j$ for $1\leq j \leq P$ and $\tilde c_n(t)=0$ otherwise). Also, let $c(t)=(c_n(t))$ be the solution to $(\mathcal{F}NLS)$ with initial data $c(0)=A+\tilde c(0)$ (i.e. $c_0(0)=A, c_n(0)=\tilde c_n(0)$ for $n \neq 0$). Since $\|b(t)\|_{l^\infty}\sim 1$, we have that $\|\tilde c(t)\|_{l^\infty} \sim \lambda^{-1}$ for all $0\leq t \leq \lambda^2T_0$. 

It will be convenient to adopt the following norm:
\begin{equation}\label{def of h norm}
\|(a_n)\|_{h^s(\Z^2)}:= \left (\sum_{n \in \Z^2}\langle n \rangle^{2s}|a_n|^2\right)^{1/2}\sim \|\sum_{n \in \Z^2} a_n e^{in.x}\|_{H^s(\T^2)}.
\end{equation}
We will often abuse notation and refer to it as the Sobolev norm of $(a_n)$.

The idea is to verify that 1) the Sobolev norm of $\tilde c(t)$ grows from $\delta$ to $K$ in the time interval $[0,\lambda^2 T_0]$ and 2) the distance in $H^s$ between $c|_\Lambda(t)$(i.e. restricting to the frequencies of $c$ in $\Lambda$) and $\tilde c(t)$ remains $\leq 1$ over this time interval. The latter is done using the stability lemma \ref{stability lemma}. 

In the notation of Lemma \ref{stability lemma}, $B$ will be taken to be $C(P)\lambda$ (since $\Lambda$ has $C(P)$ elements and $\tilde c(t)\lesssim \lambda^{-1}$) and $\sigma$ chosen as mentioned above so that $s(1+\sigma)<1$. The condition that $1\ll B^{1+\sigma} \ll N$ will be guaranteed because $\lambda$ will be chosen below to be  $\sim_{\delta, K,P} N^s$ (and $N^{s(1+\sigma)} \ll N$). Also, with $\lambda$ large enough, one also guarantees that the stability interval $B^2\log B$ from the stability lemma \ref{stability lemma} is longer than $\lambda^2 T_0\sim_{\delta, K} B^2$, which is the time needed for the Sobolev norm of $(c_n(t))$ to grow from $\delta$ to $K$.

We start by verifying that the Sobolev norm of $\tilde c(t)=(\tilde c_n(t))$ grows from $\delta$ to $K$ in the time interval $[0,\lambda^2 T_0]$. We will do this by first making sure that $\|\tilde c(0)\|_{h^s(\Z^2)}\sim  \delta$ (this will specify $\lambda$ in terms of $N$) and then verifying that the ratio of Sobolev norms 

\begin{equation}\label{ratio of Sobolev norms}
\frac{\|\tilde c(\lambda^2T_0)\|_{h^s(\Z^2)}}{\|\tilde c(0)\|_{h^s(\Z^2)}}\gtrsim \frac{K}{\delta}.
\end{equation}

Notice that since 
$$\|\tilde c (0)\|_{h^s(\Z^2)}\sim \left(\sum_{n \in \Lambda} |n|^{2s}|\tilde c_n(0)|^2\right)^{1/2}\sim_{P,\delta,K } N^s\lambda^{-1},
$$
one can make sure that $\|\tilde c(0)\|_{h^s(\Z^s)}\sim \delta$ by choosing $\lambda \sim_{P, \delta,K} N^s$. We remark that this implies also that $\| c(0)-A\|_{h^s(\Z^2)}\sim \delta$. Moving to $\eqref{ratio of Sobolev norms}$: Let $Q$ be the ratio of the Sobolev norm of $\tilde c$ at times $\lambda^2T_0$ to that at time $0$,

\begin{align*}
Q\sim&\frac{\sum_{n\in \Lambda}|n|^{2s}|\tilde c_n(\lambda^2T_0)|^2}{\sum_{n\in \Lambda}|n|^{2s}|\tilde c_n(0)|^2}\\
=& \frac{\sum_{j=1}^P\sum_{n\in \Lambda_j}|n|^{2s}|b^\lambda_j(\lambda^2T_0)|^2}{\sum_{j=1}^P\sum_{n\in \Lambda_j}|n|^{2s}|b^{\lambda}_j(0)|^2}
\end{align*}
where we have used the fact that $\tilde c(t)$ is supported on $\Lambda$ and that it comes from the diffusion solution $b(t)$. Adopting the notation $S_j=\sum_{n\in \Lambda_j}|n|^{2s}$ and using the  properties of $b(t)$ given by proposition \ref{diffusion result} we get that:

\begin{align*}
Q\gtrsim& \frac{S_{P-2}(1-\epsilon)}{\epsilon S_1+\epsilon S_2+(1-\epsilon) S_3+\epsilon S_4+\ldots +\epsilon S_P}\\
\gtrsim& \frac{(1-\epsilon)}{(1-\epsilon)\frac{S_3}{S_{P-2}}+O_{\delta, K}(\epsilon)}\\
\gtrsim &\frac{K^2}{\delta^2}
\end{align*}
by choosing $\epsilon$ small enough depending on $P,\delta, K$.

Finally, we verify that $\|c(\lambda^2T_0)\|_{h^s(\Z^2)}\gtrsim K$ by showing that $\| c|_\Lambda(\lambda^2T_0)\|_{h^s(\Z^2)} \gtrsim K$. This follows from \eqref{stability estimate} as follows:  

\begin{align*}
\left(\sum_{n \in \Lambda}|n|^{2s}|c_n(\lambda^2T_0)-\tilde c_n(\lambda^2 T_0)|^2\right)^{1/2}\leq& C(P)N^s\|c(\lambda^2T_0)-\tilde c(\lambda^2T_0)\|_{l^1(\Lambda)}\\
\leq& C(P,A) N^s\frac{\lambda^{1+\sigma}}{N^2}\leq C(P,A,\delta) \frac{N^{s(2+\sigma)}}{N^2}\\
\leq& 1\quad \hbox{if $N$ is large enough},
\end{align*}
where we have used the fact that $B=C(P)\lambda$ in the second inequality, that $\lambda=C(P,\delta) N^s$ in the third, and that $2>s(2+\sigma)$ in the last inequality. 

\endproof

\subsection{Proof of stability lemma \ref{stability lemma}}\label{proof of stability lemma}

Denote by $\mathcal{D}^*=\mathcal{D}\setminus 0$. We would like to show that $\left(c_n(t)\right)$ remains close to $\left(\tilde{c}_n(t)\right)$ in $l^1(\mathcal{D}^*)$ for a time interval $[0,T]$ under the assumption that $\sum_{n \neq 0}|\tilde{c}_n(t)| \lesssim B^{-1}$ for some $B \gg 1$ and $T \ll B^2 \log B$.

Denoting $w_n=c_n-\tilde{c}_n$ for $n \neq 0$ we get the following
difference equation:
\begin{equation}\label{w_n}
\begin{split}
-i\partial_t w_n(t)=&-\left(|c_n(t)|^2c_n(t)-|\tilde{c}_n(t)|^2\tilde{c}_n(t)\right)\\
&+\sum_{\Gamma_{res}(n)\cap \Lambda^3}\left( c_{n_1}\overline{c_{n_2}(t)}c_{n_3}(t)-\tilde{c}_{n_1}\overline{\tilde{c}_{n_2}(t)}\tilde{c}_{n_3}(t)\right)\\
&+c_0(t)^2\overline{c_{-n}(t)}e^{-2i|n|^2t}+2\left(\sum_{n_2\neq 0,-n}\overline{c_{n_2}}c_{n+n_2}e^{2i n\cdot n_2t}\right)c_0(t)\\
&+\left(\sum_{n_1\neq 0,n}c_{n_1}(t)c_{n-n_1}(t)e^{2in_1\cdot (n_1-n)t}\right)\overline{c_0}(t)\\
&+\sum_{\tilde{\Gamma}(n)\setminus \left(\Gamma_{res}(n)\cap \Lambda^3\right)}c_{n_1}\overline{c_{n_2}}c_{n_3}e^{i\omega_4t}
\end{split}
\end{equation}

Defining:

\begin{align*}
A_n(t)=&|c_n(t)|^2c_n(t)-|\tilde{c}_n(t)|^2\tilde{c}_n(t), \quad B_n(t)=\sum_{\Gamma_{res}(n)\cap \Lambda^3}\left( c_{n_1}\overline{c_{n_2}(t)}c_{n_3}(t)-\tilde{c}_{n_1}\overline{\tilde{c}_{n_2}(t)}\tilde{c}_{n_3}(t)\right)\\
C_n(t)=&c_0(t)^2\overline{c_{-n}(t)}e^{-2i|n|^2t}, \quad D_n(t)=\left(\sum_{n_2\neq 0,-n}\overline{c_{n_2}}c_{n+n_2}e^{2i n\cdot n_2t}\right)c_0(t)\\
E_n(t)=&\left(\sum_{n_1\neq 0,n}c_{n_1}(t)c_{n-n_1}(t)e^{2in_1\cdot (n_1-n)t}\right)\overline{c_0}(t), \hbox{ and }
F_n(t)=\sum_{\tilde{\Gamma}(n)\setminus \left(\Gamma_{res}(n)\cap \Lambda^3\right)}c_{n_1}\overline{c_{n_2}}c_{n_3}e^{i\omega_4t}
\end{align*}

we get that the equation satisfied by $\left(w_n(t)\right)$ can be written in integral form as:

\begin{equation}\label{w_n integral compact}
w_n(t)=i\int_0^t\left(A_n(t)+B_n(t)+C_n(t)+D_n(t)+E_n(t)+F_n(t)\right)dt
\end{equation}

All the sums below are over $\mathcal{D}^*:=\D \setminus 0$ unless otherwise indicated and implicit constant are allowed to depend on $A\sim 1$. 

The proof follows a bootstrap argument: We will first adopt the following bootstrap hypotheses: 
\begin{equation}\label{bootstrap hypothesis}
\sum_{n \neq 0}|w_n(t)|\lesssim B^{-1}
\end{equation}
for all $0 \leq t \leq T$ (and hence $\sum_{n\in \Lambda}|c_n(t)|\lesssim B^{-1}$ over the same time interval). The goal will then be to prove that $\sum_{n \neq 0}|w_n(t)|$ actually satisfies the stronger bound 

\begin{equation}\label{bootstrap conclusion}
\sum_{n \neq 0} |w_n(t)| \lesssim_A \frac{B^{1+\sigma}}{N^2} \ll B^{-1-\sigma}
\end{equation}
on $[0,T]$ if $T \ll B^2 \log B$.

Before we start estimating the contribution of the terms $A_n, \ldots, F_n$, we will single out a useful identity that we will use repeatedly:
\begin{lemma}
For any $0\leq t \leq T$,

\begin{equation}\label{c_0 estimate}
\left|A^2-|c_0(t)|^2\right| \lesssim B^{-2}
\end{equation}
In particular, $|c_0(t)|\sim A\sim_A 1$ for all $0\leq t \leq T$.
\end{lemma}
\proof
This follows directly from the bootstrap hypothesis and the conservation of mass $M[u(t)]=\sum_{n \in \mathcal{D}}|c_n(t)|^2$ which gives:

$$A^2+\sum_{n\neq 0}|c_n(0)|^2=M[u(0)]=M[u(t)]=|c_0(t)|^2+\sum_{n \neq 0}|c_n(t)|^2
$$
\endproof

\subsubsection{Contribution of $A_n$:}

Recall that
$$
A_n(t)=|c_n(t)|^2c_n(t)-|\tilde{c}_n(t)|^2\tilde{c}_n(t)=|c_n(t)|^2w_n(t)+c_n(t)\overline{w_n(t)}\tilde{c}_n(t)
+|\tilde{c}_n(t)|^2w_n(t)
$$
and hence:
\begin{equation}\label{Contribution of A_n}
\sum_{n \neq 0}|\int_0^t A_n(s)ds| \lesssim \sum_{n\neq 0} B^{-2}\int_0^t|w_n(s)|ds
\end{equation}

\subsubsection{Contribution of $B_n$:}
Recall that 
$$
B_n(t)=\sum_{\Gamma_{res}(n)\cap \Lambda^3}\left( c_{n_1}\overline{c_{n_2}(t)}c_{n_3}(t)-\tilde{c}_{n_1}\overline{\tilde{c}_{n_2}(t)}\tilde{c}_{n_3}(t)\right)
$$
Hence,
\begin{align*}
B_n(t)=&\sum_{\Gamma_{res}(n)\cap \Lambda^3}\left( c_{n_1}\overline{c_{n_2}(t)}c_{n_3}(t)-\tilde{c}_{n_1}\overline{\tilde{c}_{n_2}(t)}\tilde{c}_{n_3}(t)\right)\\
&=\sum_{\Gamma_{res}(n)\cap \Lambda^3}\left( w_{n_1}\overline{c_{n_2}(t)}c_{n_3}(t)+c_{n_1}\overline{w_{n_2}(t)}\tilde{c}_{n_3}(t)+\tilde{c}_{n_1}\overline{\tilde{c}_{n_2}(t)}w_{n_3}(t)\right)\\
\end{align*}
and we also get:

\begin{equation}\label{Contribution of B_n}
\sum_{n \neq 0}|\int_0^t B_n(s)ds| \lesssim \sum_{n\neq 0} B^{-2}\int_0^t|w_n(s)|ds
\end{equation}
\subsubsection{Contribution of $C_n$:}
Recall that
$$
C_n(t)=c_0(t)^2\overline{c_{-n}(t)}e^{-2i|n|^2t}
$$
which gives upon integration by parts:
\begin{align}
\label{bndry terms} \int_0^t C_n(s)ds =& \frac{-1}{2i|n|^2}c_0(s)^2\overline{c_{-n}(s)}e^{-2i|n|^2s}\bigg|_0^t\\
+& \frac{1}{2i|n|^2}\int_0^t\left(2c_0(s)\partial_tc_0(s)\overline{c_{-n}(s)}
+c_0(s)^2\partial_t\overline{c_{-n}(s)}\right)e^{-2i|n|^2s}ds
\end{align}

Note that since $n\in \D \setminus 0$, we have $|n|\geq N$ and hence the boundary terms \eqref{bndry terms} have the acceptable contribution of $O(B^{-1}N^{-2})$. Denote by:
\begin{align*}
I_n(t)=&\int_0^t2c_0(s)\partial_tc_0(s)\overline{c_{-n}(s)}
e^{-2i|n|^2s}ds, \quad \hbox{and} \quad II_n(t)=&\int_0^tc_0(s)^2\partial_t\overline{c_{-n}(s)}e^{-2i|n|^2s}ds.
\end{align*}
We start by estimating $I_n(t)$. From $\eqref{c_0}$, we have
\begin{equation}\label{crude c_0 getaround}
\begin{split}
\partial_t c_0(t)=&-i|c_0|^2c_0(t) +O(AB^{-2})=-iA^2c_0(t)+i(A^2-|c_0(t)|^2)c_0(t)+O(B^{-2})\\
=&-iA^2c_0(t)+O(AB^{-2})
\end{split}
\end{equation}
where we used $\eqref{c_0 estimate}$ in the last equality. This leads to:
\begin{align*}
|I_n(t)|\sim&|\int_0^tc_0(s)\partial_tc_0(s)\overline{c_{-n}(s)}
e^{-2i|n|^2s}ds|\\
\lesssim& |\int_0^tc_0(s)^2\overline{c_{-n}(s)}e^{-2i|n|^2s}ds|+B^{-2}\int_0^t |c_{-n}(s)|ds\\
=&|\int_0^t C_n(s) ds|+B^{-2}\int_0^t |c_{-n}(s)|ds, 
\end{align*}
which gives that:
\begin{equation}\label{I_n estimate}
\sum_{n\in \D^*}\frac{1}{|n|^2}|I_n(t)|\lesssim \frac{1}{N^2} \sum_{n \in \D^*} |\int_0^t C_n(s) ds|+ \frac{B^{-2}}{N^2}\sum_{n \in \D^*}\int_0^t |c_{-n}(s)|ds
\end{equation}

To estimate $II_n(t)$, we first use the expression of $\partial_t \overline{c_{-n}}$ from \eqref{c_n} to get:
\begin{align}\label{II_n expansion}
II_n(t)=& \int_0^tc_0(s)^2\partial_t\overline{c_{-n}(s)}e^{-2i|n|^2s}ds\\
\label{Xi_1}=& i\int_0^tc_0(s)^2|c_{-n}(s)|^2\overline{c_{-n}(s)}e^{-2i|n|^2s}ds\\
\label{Xi_2}&-i\int_0^t|c_0(s)|^4c_{n}(s)ds\\
\label{Xi_3}&-2i\int_0^t|c_0(s)|^2c_0(s)\left(\sum_{n_2\neq 0,n}c_{n_2}\overline{c_{-n+n_2}}e^{2i n\cdot n_2s}\right)e^{-2i|n|^2s}ds\\
\label{Xi_4}&-i\int_0^tc_0(s)^3\left(\sum_{n_1\neq 0,-n}\overline{c_{n_1}(s)c_{-n-n_1}(s)}e^{-2in_1\cdot (n_1+n)s}\right)e^{-2i|n|^2s}ds\\
\label{Xi_5}&-i\int_0^tc_0(s)^2\sum_{\tilde{\Gamma}(-n)}\overline{c_{n_1}}c_{n_2}\overline{c_{n_3}}e^{-i\omega_4s}
e^{-2i|n|^2s}ds\\
=& \Phi_1+\Phi_2+\Phi_3+\Phi_4+\Phi_5
\end{align}
where $\Phi_1,\ldots,\Phi_5$ are defined by $\eqref{Xi_1}, \ldots, \eqref{Xi_5}$ respectively.

$\sum_{n\neq 0}|\Phi_1|$ and $\sum_{n\neq 0}|\Phi_5|$ can be easily estimated by $\lesssim_A B^{-2}\sum_{n \in \D^*}\int_0^t|c_n(s)|ds$. Also, $\sum_{n\neq 0}\Phi_2 \lesssim_A \sum_{n \in \D^*} \int_0^t|c_n(s)|ds$ and 
$$\left|\sum_{n\neq 0} \Phi_3\right|, \left|\sum_{n\neq 0} \Phi_4\right|\lesssim B^{-1}\sum_{n\neq 0} \int_0^t|c_n(s)|ds.
$$ 
As a result, we get:
\begin{equation}\label{II_n estimate}
\sum_{n\neq0}|II_n(t)|\lesssim \int_0^t\sum_{n\neq0}|c_n(s)|\,ds
\end{equation}

\begin{remark}The term \eqref{Xi_2} gives the worst contribution in all our estimates and it is the reason why the result is restricted to the range $0<s<1$. It is the result of a chain of two non-resonant interactions of the zero mode with the modes at $n$ and $-n$ respectively (via the term $c_0^2 \overline{c_{n}}e^{-2i|n|^2t}$ that shows up in the equation of $-i\partial_t c_{-n}$, and the term $c_0^2 \overline{c_{-n}}e^{-2i|n|^2t}$ in the equation of $-i\partial_t c_{n}$). The net effect of this interaction chain is a \emph{second-order resonant} self-excitation of the mode at $n$.
\end{remark}

We are now in a good position to sum the contribution of $\int_0^tC_n(t)$. Using $\eqref{bndry terms}, \eqref{I_n estimate},$ and $\eqref{II_n estimate}$, and absorbing the first term of \eqref{I_n estimate} in the left-hand side we get:
\begin{equation}\label{Contribution of C_n}
\begin{split}
\sum_{n\neq 0}\bigg|\int_0^t C_n(t)dt\bigg|\lesssim&\frac{1}{N^2}\int_0^t\sum_{n\neq 0} |c_n(s)|ds
+\frac{B^{-1}}{N^2} 
\lesssim \frac{1}{N^2}\int_0^t\sum_{n\neq 0} |w_n(s)|ds
+\frac{B^{-1}+B^{-1}t}{N^2}
\end{split}
\end{equation}

\subsubsection{Contribution of $D_n$:} $D_n$ involves resonant and
non-resonant interactions, which we will separate as follows:
\begin{equation}\label{D_n=tilde+underline}
\begin{split}
D_n(t)=&\left(\sum_{n_2\neq 0,-n}\overline{c_{n_2}}c_{n+n_2}e^{2i n\cdot n_2t}\right)c_0(t)\\
&=\left(\sum_{\substack{n_2\neq 0, -n\\n_2\cdot n=0}}\overline{c_{n_2}}c_{n+n_2}\right)c_0(t)+
\left(\sum_{\substack{n_2\neq 0, -n\\n_2\cdot n\neq0}}\overline{c_{n_2}}c_{n+n_2}e^{2i n\cdot n_2t}\right)c_0(t)\\
&=\tilde{D}_n(t)+\underline{D}_n(t)
\end{split}
\end{equation}
We start by estimating $\tilde{D}_n(t)$. For this we write first:
\begin{align*}
c_0(t)=&\frac{1}{A^2}A^2c_0(t)=\frac{1}{A^2}|c_0(t)|^2c_0(t)+(A^2-|c_0(t)|^2)c_0(t)\\
=&\frac{1}{A^2}(i\partial_tc_0(t))+O(AB^{-2})+(A^2-|c_0(t)|^2)c_0(t)\\
=&\frac{1}{A^2}(i\partial_tc_0(t))+O(AB^{-2})
\end{align*}
where we have used $\eqref{c_0}$ and $\eqref{c_0 estimate}$. As a result, we get:
\begin{align*}
\int_0^t \tilde{D}_n(s)ds=&\sum_{\substack{n_2\neq 0, -n\\n_2\cdot n=0}}\int_0^t \overline{c_{n_2}(a)}c_{n+n_2}(a)\left(\frac{1}{A^2}(i\partial_tc_0(s))+O(B^{-2})\right)ds\\
=&\frac{i}{A^2}\sum_{\substack{n_2\neq 0, -n\\n_2\cdot n=0}}\bigg(\overline{c_{n_2}(.)}c_{n+n_2}(.)c_0(.)\bigg|_0^t\\
&-\int_0^t\left(\partial_t\overline{c_{n_2}(s)}c_{n+n_2}(s)c_0(s)+\overline{c_{n_2}(s)}\partial_tc_{n+n_2}(s)c_0(s)\right)ds\bigg)\\
&+O\left(B^{-2}\sum_{\substack{n_2\neq 0, -n\\n_2\cdot n=0}}\int_0^t |c_{n_2}(s)||c_{n+n_2}(s)|ds\right)
\end{align*}

The contribution of the boundary terms above can be readily estimated in $l^1(\D^*)$ by $\lesssim B^{-2}$. For the first integral, we again use equation $\eqref{c_n}$ for $\partial_t c_{n_2}$ and $\partial_t c_{n-n_2}$ to get (we do the calculation for $\partial_tc_{n_2}$, that of $\partial_tc_{n-n_2}$ is similar):

\begin{align*}
&\int_0^t\partial_t\overline{c_{n_2}(s)}c_{n+n_2}(s)c_0(s)ds\sim \int_0^t|c_0(s)|^2\overline{c_0(s)}c_{-n_2}(s)c_{n+n_2}(s)e^{2i|n_2|^2s}ds\\
&-\int_0^t c_0(s)c_{n+n_2}(s)|c_{n_2}|^2\overline{c_{n_2}} ds +2\int_0^t |c_0(s)|^2c_{n+n_2}(s)\left(\sum_{l\neq 0,-{n_2}}c_{l}\overline{c_{{n_2}+l}}e^{-2i {n_2}\cdot ls}\right) ds\\
&+\int_0^t c_0(s)^2c_{n+n_2}(s)\left(\sum_{l\neq 0,{n_2}}\overline{c_{l}}\overline{c_{{n_2}-l}}e^{-2il\cdot (l-{n_2})s}\right)ds+\int_0^t c_0(s)c_{n+n_2}(s)\left(\sum_{\tilde{\Gamma}({n_2})}\overline{c_{l_1}}c_{l_2}\overline{c_{l_3}}e^{-i\omega_4s}\right)ds\\
=& \Xi_1 +\Xi_2+\ldots+\Xi_5
\end{align*}
where 
\begin{align*}
&\Xi_1=\int_0^t|c_0(s)|^2\overline{c_0(s)}c_{-n_2}(s)c_{n+n_2}(s)e^{2i|n_2|^2s}ds, \quad \Xi_2=-\int_0^t c_0(s)c_{n+n_2}(s)|c_{n_2}|^2\overline{c_{n_2}} ds\\
&\Xi_3=2\int_0^t c_0(s)c_{n+n_2}(s)\left(\sum_{l\neq 0,-{n_2}}c_{l}\overline{c_{{n_2}+l}}e^{-2i {n_2}\cdot ls}\right)\overline{c_0(s)}ds, \\
&\Xi_4= \int_0^t c_0(s)c_{n+n_2}(s)\left(\sum_{l\neq 0,{n_2}}\overline{c_{l}}\overline{c_{{n_2}-l}}e^{-2il\cdot (l-{n_2})s}\right)c_0 ds\\
&\Xi_5= \int_0^t c_0(s)c_{n+n_2}(s)\bigg(\sum_{\tilde{\Gamma}({n_2})}\overline{c_{l_1}}c_{l_2}\overline{c_{l_3}}e^{-i\omega_4s}\bigg)ds
\end{align*}

By Property VIII$_\Lambda$, either $n_2$ or $n+n_2 \notin \Lambda$: if $n+n_2 \notin \Lambda$, then $c_{n+n_2}=w_{n+n_2}$ and we easily estimate:

\begin{equation}\label{case n+n_2 notin Lambda}
\left|\sum_{\substack{n,n_2 \neq 0\\ n.n_2=0}} \int_0^t\partial_t\overline{c_{n_2}(s)}c_{n+n_2}(s)c_0(s)ds\right| \lesssim B^{-2}\int_0^t\sum_{n\neq 0}|w_n(s)|ds
\end{equation}
where we integrated by parts in $\Xi_1$ (to gain a factor of $N^{-2}\ll B^{-2}$) and used direct estimates for $\Xi_2, \ldots, \Xi_5$.

The case when $n+n_2\in \Lambda$ ($\Rightarrow n_2 \notin \Lambda$) is slightly more tedious to bound and we need to consider the different contributions of $\Xi_1, \ldots \Xi_5$. We start by estimating $\Xi_1$: Integrating by parts and using $\eqref{c_0}$ and $\eqref{c_n}$ to estimate $\partial_t c_0=O(1)$ and $\sum_{n \in \D^*}|\partial_t c_n| \lesssim B^{-1}$, we arrive at the estimate:

\begin{equation}\label{estimate on Xi_1}
\sum_{\substack{n,n_2 \neq 0\\ n.n_2=0}}|\Xi_1|\lesssim \frac{B^{-1}}{N^2}\int_0^t \sum_{n \neq 0} |c_n(s)|ds +\frac{B^{-2}}{N^2}\lesssim \frac{B^{-1}}{N^2}\int_0^t \sum_{n \neq 0} |w_n(s)|ds +\frac{B^{-2}(1+t)}{N^2}.
\end{equation}

Obviously, since $n_2 \notin \Lambda$, we readily have the estimate:
\begin{equation}\label{estimate on Xi_2}
\sum_{\substack{n,n_2 \neq 0\\ n.n_2=0}}|\Xi_2|\lesssim B^{-3}\int_0^t \sum_{n\neq 0}|w_n(s)|ds
\end{equation}

As for $\Xi_3$, we split the sum over $l$ into the resonant contribution when $l.n_2=0$ and the non-resonant contribution when $l.n_2 \neq 0$. For the resonant part, we use Property VIII$_\Lambda$ to conclude that either $l$ or $n_2+l \notin \Lambda$ and hence we estimate the contribution of this part by:

$$
\lesssim B^{-2}\int_0^t \sum_{n \neq 0}|w_n(s)|ds
$$

As for the non-resonant part of $\Xi_3$ when $n_2.l \neq 0$, we integrate by parts to gain a factor of $n_2.l$ in the denominator (recall that since $n_2, l \in \mathcal{D}, |n_2. l|\geq N^2$) and estimate the resulting contribution using $\eqref{c_0}$ and $\eqref{c_n}$ by: 
\begin{equation}\label{estimate on Xi_3}
\begin{split}
\sum_{\substack{n,n_2 \neq 0\\ n.n_2=0}}|\Xi_3|\lesssim& \frac{B^{-3}}{N^2}+\frac{B^{-2}}{N^2}\int_0^t \sum_n |c_n(s)| ds\\
\lesssim& \frac{B^{-2}}{N^2}\int_0^t \sum_{n\neq 0}|w_n(s)|ds+\frac{B^{-3}(1+t)}{N^2}
\end{split}
\end{equation}

The estimate on $\Xi_4$ is similar to that on $\Xi_3$ and we won't repeat the calculation. Finally, to estimate the contribution of $\Xi_5$, we notice that if one of $l_1, l_2, l_3 \notin \Lambda$, we directly have the estimate:
$$
\sum_{\substack{n,n_2 \neq 0\\ n.n_2=0}}|\Xi_5|\lesssim B^{-3}\int_0^t \sum_{n\neq 0}|w_n(s)|ds
$$ 
Otherwise, $l_1, l_2, l_3 \in \Lambda$ and as a result of Property II$_\Lambda$ and the fact that $n_2 \notin \Lambda$, we get that $w_4 \neq 0$ which implies that $|\omega_4|>N^2$. Integrating by parts and estimating as before we arrive at the estimate:
\begin{equation}\label{estimate on Xi_5}
\sum_{\substack{n,n_2 \neq 0\\ n.n_2=0}}|\Xi_5|\lesssim \frac{B^{-3}}{N^2}\int_0^t \sum_{n\neq 0}|w_n(s)|ds+\frac{B^{-4}(1+t)}{N^2}
\end{equation}

Combining the estimates on $\Xi_1 \ldots \Xi_5$ in $\eqref{case n+n_2 notin Lambda}-\eqref{estimate on Xi_5}$ we get the following estimate on $\int_0^t \tilde{D}_n(t)dt$:
\begin{equation}\label{tilde D_n estimate}
\sum_{n\neq 0}|\int_0^t \tilde{D}_n(s)ds|\lesssim  B^{-2}\int_0^t \sum_{n\neq 0}|w_n(s)|ds+\frac{B^{-2}(1+t)}{N^2}
\end{equation}
where we used the fact that $B \lesssim N$.

To finish estimating the contribution of $D_n$, we are left with its non-resonant part
$\underline{D}_n(t)$ of $D_n(t)$ given by:
$$
\underline{D}_n(t)=\sum_{n_2\cdot n \neq 0, n_2 \neq -n}\overline{c_{n_2}}c_{n+n_2}e^{2i n\cdot n_2t}c_0(t)\\
$$
Integrating by parts and estimating as before we get:

\begin{equation}\label{underline D_n estimate}
\sum_{n\neq 0}|\int_0^t \underline{D}_n(s)dt|\lesssim  \frac{B^{-2}}{N^2}+\frac{B^{-1}}{N^2}\sum_{n\neq 0}\int_0^t |c_n(s)|ds
\end{equation}
and consequently, we finally obtain:
\begin{equation}\label{Contribution of D_n}
\sum_{n\neq 0}|\int_0^t D_n(s)ds|\lesssim  B^{-2}\int_0^t \sum_{n\neq 0}|w_n(s)|ds+\frac{B^{-2}(1+t)}{N^2}
\end{equation}

\subsubsection{Contribution of $E_n(t)$:} Estimating the contribution of $E_n(t)$ is similar to that of $D_n(t)$ and we won't repeat the argument that gives:
\begin{equation}\label{Contribution of E_n}
\sum_{n\neq 0}|\int_0^t E_n(s)ds|\lesssim  B^{-2}\int_0^t \sum_{n\neq 0}|w_n(s)|ds+\frac{B^{-2}(1+t)}{N^2}
\end{equation}

\subsubsection{Contribution of $F_n(t)$:}
Recall that
$$
F_n(t)=\sum_{\tilde{\Gamma}(n)\setminus \left(\Gamma_{res}(n)\cap \Lambda^3\right)}c_{n_1}\overline{c_{n_2}}c_{n_3}e^{i\omega_4t}
$$
Here the main concern are the terms for which $n_1,n_2,n_3 \in \Lambda$. To deal with those we integrate by parts in the integral of $F_n$ and use the fact that when $w_4\neq0$ then $|w_4|\geq N^2$. Noting also that $\|\partial_t c_{n_k}\|_{l^1(D^*)}=O(B^{-1})$ and using the bootstrap assumption we get:
\begin{align*}
\sum_{n\neq 0} |\int_0^t \left(\sum_{\tilde{\Gamma}(n)\cap \Lambda^3\setminus \left(\Gamma_{res}(n)\cap \Lambda^3\right), }c_{n_1}\overline{c_{n_2}}c_{n_3}e^{i\omega_4s}\right)ds| \lesssim& \frac{B^{-3}(1+t)}{N^2}+\frac{B^{-2}}{N^2}\int_0^t\sum_{n\neq 0}|w_n(s)|ds
\end{align*}

As for the other terms of $F_n(t)$, we will always have at least one of $n_1,n_2,n_3$ not in $\Lambda$ and hence this sum can be estimated directly by $B^{-2}\int_0^t|w_n(s)|ds$. Combining this with the above yields the estimate:

\begin{equation}\label{Contribution of F_n}
\sum_{n\neq 0} |\int_0^t F_n(s)ds| \lesssim \frac{B^{-3}(1+t)}{N^2}+B^{-2}\sum_{n\neq 0}\int_0^t|w_n(s)|ds
\end{equation}

\subsubsection{End of proof:}

Summing the contributions of $A_n ,\ldots, F_n$ in \eqref{Contribution of A_n}, \eqref{Contribution of B_n}, \eqref{Contribution of C_n}, \eqref{Contribution of D_n}, \eqref{Contribution of E_n}, \eqref{Contribution of F_n}, we get that:
$$
\sum_{n \neq 0}|w_n(s)|ds \lesssim B^{-2}\int_0^t \sum_{n \neq 0}|w_n(s)|ds+\frac{B^{-1}(1+t)}{N^2}
$$
where we used the fact that $B\lesssim N$. Gronwall's inequality now gives:
\begin{equation}\label{|w_n| for small s}
\sum_{n\neq 0}|w_n|\lesssim (\frac{B^{-1}(1+t)}{N^2})e^{CB^{-2}t}
\end{equation}
and so for $0 \leq t \leq T \ll_\sigma B^2\log B$, we have:
\begin{equation}\label{|w_n| for small s}
\sum_{n\neq 0}|w_n|\lesssim (\frac{B^{-1}B^{\sigma/2}(1+t)}{N^2})\ll \frac{B^{1+\sigma}}{N^2}\ll B^{-1-\sigma}
\end{equation}
for $B$ large enough and $N\gg B^{1+\sigma}$.
\endproof

\section{Unbounded orbits for almost-polynomial nonlinearities}\label{LT instability section}

In this section, we consider the equation:
\begin{equation}\label{NLSR}
\begin{split}
-i\partial_t u +\Delta u=&\mathcal{N}_R(u)\\
u(0,x)=& u_0(x)
\end{split}
\end{equation}
with
\begin{equation}\label{def of N_R}
\mathcal{N}_R(u)=\sum_{n \in \Z^d} \left(\sum_{\Gamma_R(n)} \widehat{u}(n_1) \overline{\widehat{u}(n_2)} \widehat u(n_3)\right) e^{2\pi i n.x}
\end{equation}
and $\Gamma_R(n)=\{ (n_1, n_2, n_3) \in (\Z^d)^3: n_1-n_2+n_3=n \hbox{ and } \omega_4:= |n_1|^2-|n_2|^2+|n_3|^2-|n|^2 \in [-R, R]\}$. This is obtained from the cubic nonlinearity (which corresponds to the case $R=\infty$) by deleting highly non-resonant interactions for which $|\omega_4|> R$.

We will prove that for any $R<\infty$, \eqref{NLSR} admits unbounded orbits by following the program suggested by Proposition \ref{observation}. Our results will be valid for all Sobolev spaces $H^s(\T^2)$ with $s>0, s\neq 1$\footnote{In the range $0<s<1$, the result will be conditional on small-data global well-posedness, which can be proved in some cases at least when the mass $M$ is below a certain threshold  (cf. following remark).}. 

We start by making some remarks about the nonlinearity $\mathcal{N}_R(u)$. We first note that for any $R$,  one can easily verify that the flow defined by \eqref{NLS} conserves the mass:
$$
M[u(t)]:= \int_{\T^2} |u(t,x)|^2 dx.
$$
Another conserved quantity is the Hamiltonian:
\begin{equation}\label{Hamiltonian}
H_R(u)= \sum_{n}|n|^2 |\widehat u(n)|^2 + \frac{1}{2}\sum_{\substack{n_1-n_2+n_3-n_4=0\\ \left||n_1|^2-|n_2|^2+|n_3|^2-|n_4|^2\right| \leq R}} \widehat{u}(n_1)\overline{\widehat{u}(n_2)}\widehat{u}(n_3) \overline{\widehat{u}(n_4)}.
\end{equation}

Notice that when $R=0$, this is the Hamiltonian of the resonant system (or $(\mathcal{RF}NLS)$) from the previous sections, and in this case $H_{R=0}(u)$ becomes coercive, in the sense that it controls the $\dot H^1$ norm of the solution. Indeed, the second term in \eqref{Hamiltonian} can be written as:
\begin{equation}\label{R=0 Hamiltonian}
\sum_{n\in \Z^2, k \in \Z} \left| \sum_{\substack{n_1+n_3=n\\ |n_1|^2+|n_3|^2=k}} \widehat u(n_1) \widehat u(n_3)\right|^2\geq 0
\end{equation}

\begin{remark}
(\emph{Coercivity of $H_R$ and global existence}) The local well-posedness theory of \eqref{NLSR} is exactly the same as that of the polynomial cubic nonlinearity $|u|^2u$ in \eqref{cubic NLS}. This was proved in \cite{Thesis} (chapter 10.1) by adapting Bourgain's proof for the polynomial cubic nonlinearity in \cite{B1}. Indeed, just like the cubic nonlinearity $|u|^2u$, \eqref{NLSR} is locally well-posed in $H^s(\T^2)$ for all $s>0$ and the time of existence depends only on $\|u_0\|_{H^s(\T^2)}$ for $0<s<1$ and on $\|u_0\|_{H^1(\T^2)}$ for $s\geq 1$. 

When $R=0$ or $\infty$, one can directly combine this local well-posedness result with the a priori control on the $H^1$ norm (given by the conservation and coercivity of the Hamiltonian $H_R$) to obtain global well-posedness of \eqref{NLSR} for all $s\geq 1$. When $R\neq 0$, one can still extend the local well-posedness in $H^s(\T^2)$ for $s\geq 1$ into a global existence statement at least when the mass $M$ is smaller than some universal constant $C$. In fact, we have:
\begin{align*}
&\frac{1}{2}\sum_{\substack{n_1-n_2+n_3-n_4=0\\ \left||n_1|^2-|n_2|^2+|n_3|^2-|n_4|^2\right| \leq R}} \widehat{u}(n_1)\overline{\widehat{u}(n_2)}\widehat{u}(n_3) \overline{\widehat{u}(n_4)}\leq\frac{1}{2}\sum_{n_1-n_2+n_3-n_4=0} |\widehat{u}(n_1)|\,|\overline{\widehat{u}(n_2)}|\,|\widehat{u}(n_3)\,| \overline{\widehat{u}(n_4)}|\\
&\leq \frac{1}{2}C_{GN} \left(\sum_{n} |\widehat u (n)|^2\right)\left( \sum_{n}|n|^2 |\widehat u(n)|^2\right)
\end{align*}
where we applied the Gagliardo-Nirenberg inequality in the last step and denoted by $C_{GN}$ the optimal constant of this inequality. So if $M[u(0)] <(\frac{2}{C_{GN}})^{1/2}$, the Hamiltonian $H_R(u)$ becomes coercive (uniformly in $R$), which gives (along with the conservation of mass) the needed a priori control on the $H^1$ norm. For concreteness, we summarize the contents of this remark in the following lemma, whose full proof can be found in \cite{Thesis} (chapter 10.1).

\begin{lemma}\label{WP for NLSR}
\begin{enumerate}
\item For any $R \in {0}\cup \N \cup \{+\infty\}$, equation $\eqref{NLSR}$ is locally well-posed in $H^s(\T^2)$ for any $s>0$. Moreover, the length of the local existence interval depends on $\|u_0\|_{H^s(\T^2)}$ when $0<s<1$ and only on $\|u_0\|_{H^1(\T^2)}$ when $s\geq 1$.
\item When $R=0$ or $+\infty$, the equation is globally well-posed in $H^s(\T^2)$ for all $s\geq 1$.
\item There exists a universal constant $C^*:=(\frac{2}{C_{GN}})^{1/2}$, such that if $M[u_0]\leq C^*$, then the local solution in part i) exists globally in $H^s(\T^2)$ for any $s\geq 1$.
\end{enumerate}
\end{lemma}

Finally, we end the remark by noting that minor modifications to the arguments in \cite{Hani1, Hani2} (where the case $R=\infty$ is considered) can extend the above global well-posedness statements in parts ii) and iii) above to the range $s>2/3$.\endproof
\end{remark}
Using the same reductions as in Section \ref{reductions subsection}, we get that equation $\eqref{NLSR}$ is equivalent to the following system of ODE:
\begin{equation}\label{NNLS2}
-i\partial_t a_n=-a_n|a_n|^2+\sum_{n_1,n_2,n_3 \in \Gamma^R(n)}a_{n_1}\overline{a_{n_2}}a_{n_3}e^{i\omega_4 t}
\end{equation}
where $\Gamma^R(n)=\{(n_1,n_2,n_3)\in (\Z^2)^3: n_1-n_2+n_3=n \text{ and } n_1,n_3 \neq n: |\omega_4|\leq R\}$. Recall that at any time $t$, one can recover the solution $u$ of $\eqref{NLSR}$ from $v(t,x)=\sum_{n\in \Z^2} a_n e^{i(nx+|n^2|t)}$, and vice-versa, and the two solutions have the same Sobolev norms.

\subsection{Proof of Theorem \ref{WT for special nonlinearities}}
The non-self interactions that appear in system \eqref{NNLS2} happen in parallelograms with vertices $n,n_1,n_2, n_3$ satisfying $n_1, n_3 \neq n$ (which guarantees that the parallelogram does not degenerate to a point) and $|\omega_4|\leq R$. In what follows it will often be convenient to specify the arguments of the function $\omega_4$, so we precise:
\begin{equation}\label{def of omega_4} 
\omega_4(n_1, n_2, n_3, n_4):=|n_1|^2-|n_2|^2+|n_3|^2-|n_4|^2,
\end{equation}
and remark that when $n_4=n_1-n_2+n_3$, $\omega_4$ can be understood as a measure of an ``eccentricity" of the parallelogram with vertices $n_1, n_2, n_3, n_4$ (in the sense that the smaller $|\omega_4|$ is, the closer is the parallelogram to being a rectangle).

Our initial data to \eqref{NNLS2} will be chosen to be supported on a set $S$ satisfying the following closure property:

\textbf{$R-$closure condition:} \emph{If $n_1,n_2,n_3 \in S$ are three vertices of a parallelogram with $n_4:=n_1-n_2+n_3$ being the forth vertex, and if $\omega_4=|n_1|^2-|n_2|^2+|n_3|^2-|n_4|^2\in [-R,R]$, then $n_4 \in S$.}

This condition is simply the generalization of the closure condition from Section \ref{Nuclear families}, and it guarantees that the emanating solution $a_n(t)$ remains supported on $S$ for all later times. We include this in the following lemma:

\begin{lemma}\label{Gronwall on S}
Suppose that the initial data $a_n(0) \in l^1(\Z^2)$ is supported in a set $S$ satisfying the $R-$closure condition above. If $a_n(t)\in L_{t, loc}^\infty \left(I; l^1(\Z^2)\right)$\footnote{This condition is actually redundant as it is not hard to see that \eqref{NNLS2} is locally well-posed in $l^1(\Z^2)$.} solves $\eqref{NNLS2}$ on some time interval $I\ni 0$, then $a_n(t)$ is also supported on $S$, i.e. $a_n(t)=0$ unless $n\in S$.
\end{lemma}

\proof Let $B(t)=\sum_{n\notin S}|a_n|^2$, then 

$$
B'(t)=2\Re \sum_{n \notin S}\left(-|a_n|^4+\sum_{n_1,n_2,n_3 \in \Gamma^R(n)} \overline{a_n}a_{n_1}\overline{a_{n_2}}a_{n_3}e^{i\omega_4t}\right)
$$
By the closure property of $S$ either $n_1,n_2,$ or $n_3 \notin S$ if $n \notin S$. This gives the bound:

$$
B'(t)\leq CB(t)
$$
for all $t$ in an arbitrary compact subinterval of $I$. Since $B(0)=0$, we get the result by Gronwall's inequality.\endproof

The plan is to apply Proposition \ref{observation} (ii) to conclude the existence of unbounded Sobolev orbits of \eqref{NNLS2}. For this, we will prove long-time strong instability of the flow near all initial data whose support satisfies the $R$-closure condition. This is the content of the following main proposition:

\begin{proposition}\label{LT instability c10}
Let $s>0, s\neq 1$. Suppose that $S$ is a finite subset of $\Z^2$ satisfying the $R-$closure condition and $\phi \in H^s(\T^2)$ is supported in Fourier space on $S$ (i.e. $\widehat \phi(n)=0$ unless $n\in S$). For any $\delta \ll 1$, $K \gg 1$, there exists a solution $(c_n(t))$ of $\eqref{NNLS2}$ and a time $T>0$ such that:

\begin{equation}\label{explosion2}
\begin{split}
\|\sum_{n \in \Z^2} c_n(0) e^{in.x} -\phi\|_{H^s(\T^2)} \leq \delta\quad \hbox{ and}\\
\|\sum_{n \in \Z^2}c_n(T)e^{in.x}\|_{H^s(\T^2)}\geq K.
\end{split}
\end{equation}

Moreover, $c_n(0)$ can be chosen to be supported on a finite set satisfying the $R-$closure condition.
\end{proposition}
We delay the proof of this key proposition to the next section, but let us see how it gives a proof of Theorem \ref{WT for special nonlinearities}. For $s\geq 1$, the flow map:
\begin{align*}
X:=\{u_0\in H^s: \|u_0\|_{L^2} \leq C^*\}&\to \in C_tH^s([0,T]\times \T^2) \\
u_0 &\mapsto u(t) \hbox{ solving \eqref{NLSR} with initial data } u_0,
\end{align*}
exists and gives a well-posed dynamical system on the closed set $X\subset H^s(\T^2)$. Recall that this set is invariant under the flow thanks to the conservation of mass. As mentioned before, this restriction on the mass is included to guarantee global well-posedness of the flow and is not needed when $R=0$.

Let $D\subset X$ be the subset of all functions whose Fourier transform is supported on a finite set $S$ satisfying the $R-$closure property. Define $\mathcal{F}:=\overline{D}$ to be its closure in $X$ (or equivalently in $H^s$ since $X$ is closed). The proof now proceeds exactly as that of Proposition \ref{observation} to conclude that, since $\mathcal{F}$ is a complete metric space, the set of unbounded orbits is non-empty and co-meager in $\mathcal{F}$.

\remark Since Proposition \ref{LT instability c10} holds for all $s>0, s\neq 1$, the above proof holds for $0<s<1$ as long as we can prove that the flow map is globally well-posed on some open subset of $H^s$ containing 0. As mentioned in the previous remark, this is the case for $s>2/3$ by adapting the arguments in \cite{Hani2}. 

\subsection{Proof of Proposition \ref{LT instability c10}}
The idea of the proof of Proposition \ref{LT instability c10} is as follows: Starting with an initial data $b_n(0)=\widehat \phi(n)$ supported on a finite set $S$ satisfying the $R$-closure condition, the solution $(b_n(t))$ to $\eqref{NNLS2}$ remains supported on the set $S$ for all time. One another hand, one can construct a norm exploding solution $(a_n(t))$ to $\eqref{NNLS2}$ that is supported on a \emph{disjoint} set $\Lambda$ also satisfying the $R-$closure condition, in addition to a combinatorial condition defined in terms of how it ``interacts" with $S$ (See Lemma \ref{pasting lemma}). The set $\Lambda$ will also have to satisfy a list of conditions similar to those in Section \ref{Nuclear families}, and its construction is one of the main components of the proof. The careful choice of $\Lambda$ will allow us to ``paste" the two solutions $(a_n(t))$ and $(b_n(t))$ to obtain a solution $(c_n(t))$ satisfying \eqref{explosion2}. 

We start with specifying the above-mentioned geometric/combinatorial condition on two sets $S_1$ and $S_2$ of $\Z^2$ that allows one to paste solutions supported on $S_1$ and $S_2$.

\begin{lemma}[Pasting lemma]\label{pasting lemma}
Let $S_1$ and $S_2$ be two finite disjoint subsets of $\Z^2$satisfying the $R-$closure condition. Suppose that $S_1$ and $S_2$ are \emph{not connected by any parallelogram with eccentricity $\leq R$} in the following sense: if $n_1,n_2\in S_1, n_3\in S_2$, and $n$ is the forth vertex of the parallelogram formed by $n_1, n_2,$ and $n_3$, then:
\begin{equation}\label{nonrec}
\begin{cases}
|\omega_4(n_1,n_2, n_3, n)|\geq R \quad \hbox{if }n=n_1-n_2+n_3, \\
|\omega_4(n_1, n_3, n_2, n)|\geq R \quad \hbox{if }n=n_1-n_3+n_2.
\end{cases}
\end{equation}
and vice-versa with $S_1$ and $S_2$ interchanged.

If $(a_n(t))$ and $(b_n(t))$ are solutions to $\eqref{NNLS2}$ with $(a_n(0))$ and $(b_n(0))$ supported in $S_1$ and $S_2$ respectively, then $(a_n(t)+b_n(t))$ also solves $\eqref{NNLS2}$.
\end{lemma}

\begin{remark}
The two conditions in \eqref{nonrec} corresponds to the two ways to form a parallelogram out of the three points $n_1, n_2,$ and $n_3$.
\end{remark}

\proof The proof is by direct verification. Since $a_n(t)$ and $b_n(t)$ satisfy finite dimensional systems of ODE and since $\eqref{NNLS2}$ conserves the $L^2$ norms $\sum_n|a_n(t)|^2$ and $\sum_{n}|b_n(t)|^2$, the two solutions exist globally in time. Let $c_n(t)=a_n(t)+b_n(t)$. Since $S_1\cup S_2$ clearly satisfies the $R-$closure condition, we directly get that $c_n(t)=0$ unless $n\in S_1\cup S_2$. Say $n\in S_1$, then 
\begin{align*}
-i\partial_t c_n=&-i\partial_t a_n=-a_n|a_n|^2+\sum_{n_1,n_2,n_3 \in \Gamma^R(n)}a_{n_1}\overline{a_{n_2}}a_{n_3}e^{i\omega_4t}\\
=&-c_n|c_n|^2+\sum_{n_1,n_2,n_3 \in \Gamma^R(n)}c_{n_1}\overline{c_{n_2}}c_{n_3}e^{i\omega_4t}
\end{align*}
because the last summand would vanish unless $n_1,n_2,n_3 \in S_1\cup S_2$ and in this case we must have that $n_1,n_2,n_3 \in S_1$ since $n\in S_1$.
\endproof

As a result of the above lemma, the key point is to construct a solution $(a_n(t))$ supported on frequency set $\Lambda \subset \Z^2$ so that $\Lambda$ and $S$ are not connected by parallelograms with ``eccentricity"  $|\omega_4| \leq R$ as defined in \eqref{nonrec} of the pasting lemma above. Adopting the useful notation in \eqref{def of h norm}, we will also make sure that:
\begin{equation}\label{a_n(0)}
\|(a_n(0))\|_{h^s(\Z^2)}\leq \delta
\end{equation}
and at some later time $T$:
\begin{equation}\label{a_n(T)}
\|(a_n(t))\|_{h^s(\Z^2)}\geq K.
\end{equation}

By the pasting lemma, if $(b_n(t))$ solves $\eqref{NNLS2}$ with initial data $(\widehat \phi(n))$, then $(c_n=a_n(t)+b_n(t))$ solves $\eqref{NNLS2}$ with initial data $(a_n(0)+\widehat \phi(n))$. Clearly, this solution proves Theorem \ref{LT instability c10} since
$$
\|c_n(0)-\widehat \phi(n)\|_{h^s(\Z^2)}=\|(a_n)\|_{h^s}\leq \delta
$$
and
$$
\|c_n(T) \|_{h^s(\Z^2)}\geq \|(a_n(T)\|_{h^s(\Z^2)}\geq K.
$$

\subsubsection{Construction of $\Lambda$} We start by compiling all the properties that we would like $\Lambda$ to satisfy:

\begin{itemize}
\item Property I$_\Lambda$ (Initial data): The initial data $a_n(0)$ is supported on $\Lambda$, i.e. $a_n(0)=0$ unless $n\in \Lambda$.

\item Property II$_\Lambda$ ($R-$disjointness with $S$): $\Lambda$ and $S$ are \emph{not connected by any parallelogram with eccentricity $\leq R$} in the sense that if $n_1,n_2\in S, n_3\in \Lambda$ and $n$ is the forth vertex of the parallelogram formed by $n_1, n_2,$ and $n_3$, then \eqref{nonrec} holds; and vice-versa with $\Lambda$ and $S$ interchanged.

We would also like $\Lambda$ to satisfy the $R-$closure property so we wouldn't need to deal with interactions between $\Lambda$ and the outside world. 
\item Property III$_\Lambda$ (Closure): If $n_1,n_2, n_3 \in \Lambda$ are three vertices of a parallelogram for which: 
$$
|\omega_4|=|n_1|^2-|n_2|^2+|n_3|^2-|n|^2|\leq R,
$$
then the forth vertex $n=n_1-n_2+n_3$ is also in $\Lambda$. Equivalently, if $(n_1,n_2,n_3) \in \Gamma^{R}(n)\cap \Lambda^3$, then $n\in \Lambda$.

Now we simplify the system $\eqref{NNLS2}$ by stipulating that all parallelograms in $\Lambda$ with $|\omega_4|\leq R$ are actually rectangles:

\item Property IV$_{\Lambda}$ (Rectangular structure): If $n_1,n_2,n_3,n \in \Lambda$ are four vertices of a parallelogram with $n_1+n_3=n+n_2$ and $|\omega_4|\leq R$, then this parallelogram is actually a rectangle with $\omega_4=0$. 

This Property reduces system $\eqref{NNLS2}$ for $(a_n(t))$ into the resonant system $\eqref{RFNLS}$ that we recall here for convenience:
\begin{equation}\label{RFNLS2}
-i\partial_t a_n=-a_n|a_n|^2+\sum_{n_1,n_2,n_3 \in \Gamma_{res}(n)}a_{n_1}\overline{a_{n_2}}a_{n_3}
\end{equation}
Note that this is a finite system of ODE on $\Lambda$ since $a_n(t)=0$ unless $n \in \Lambda$.
\end{itemize}

We now assume as in Section \ref{Nuclear families} that $\Lambda$ splits into a disjoint union $\Lambda=\Lambda_1 \cup \ldots \cup \Lambda_P$ for some integer $P>1$, and we call $\Lambda_j$ the \emph{$j-$th generation} of $\Lambda$. A \emph{nuclear family} in $\Lambda$ is defined as before as a rectangle $(n_1, n_2, n_3,n_4)$ where the frequencies $n_1,n_3 $ (called the ``parents") belong to the $j-$th generation $\Lambda_j$ while the ``\emph{children}" $n_2, n_4$ belong to $\Lambda_{j+1}$. Recall that if $(n_1,n_2,n_3,n_4)$ is a nuclear family then so are the permutations $(n_3,n_2,n_1,n_4)$ ,$(n_1,n_4,n_3,n_2)$, and $(n_3,n_4,n_1,n_2)$. In addition to the above properties we also impose the same rules of genealogy on the set $\Lambda$ as in Section \ref{Nuclear families}, namely:

\begin{itemize}
\item Property V$_{\Lambda}$ (Existence and uniqueness of spouse and children).
\item Property VI$_{\Lambda}$ (Existence and uniqueness of parents and siblings).
\item Property VII$_\Lambda$ (Non-degeneracy).
\item Property VIII$_\Lambda$ (Faithfulness).
\end{itemize}

These are exactly Properties III$_\Lambda$-VI$_\Lambda$ of Section \ref{Nuclear families}. Note that by the closure property IV$_\Lambda$, Proporty VIII$_\Lambda$ implies that $\Lambda$ contains no other parallelograms with $|\omega_4|\leq R$ apart from nuclear families (which are genuine rectangles). 

Assuming that such a set exists, the system $\eqref{RFNLS2}$ can now be written as:
\begin{equation}\label{new system}
-i\partial_t a_n=-a_n|a_n|^2 +2a_{n_{child-1}}a_{n_{child-2}}\overline{a_{n_{spouse}}}+2a_{n_{parent-1}}a_{n_{parent-2}}\overline{a_{n_{sibling}}}
\end{equation}
where we used the fact that the only rectangles that contain $n$ as a vertex are those coming from nuclear families where $n$ is either a parent (second term on the RHS) or a child (third term on RHS). The constant factor 2 on the RHS comes from the trivial permutations of the nuclear families.

Finally, we reduce the resonant system furthermore by imposing intragenerational equality:
\begin{itemize}
\item Property IX$_\Lambda$ (Intragenerational equality): The function $n\mapsto a_n(0)$ is constant on each generation $\Lambda_j$. In other words, $a_{n_1}(0)=a_{n_2}(0)$ if $n_1, n_2 \in \Lambda_j$.
\end{itemize}

One can verify (either by a Gronwall or a bootstrap argument) that Property IX$_\Lambda$ is preserved by the flow of $\eqref{RFNLS2}$ in the sense that intragenerational equality remains true for all time if it is true initially. The reason for this is that every mode in generation $\Lambda_j$ is excited in exactly the same way as other modes in generation $\Lambda_j$ (cf. the RHS of \eqref{new system}). Property IX$_\Lambda$ allows to collapse system $\eqref{RFNLS2}$, which is a system in $(a_n(t))$ with $n\in \Lambda=\Lambda_1 \cup \ldots \cup \Lambda_P$, into a smaller \emph{``Toy system"} of ODE in $b_j(t)$ with $1\leq j \leq N$ by setting $a_n(t)=b_j(t)$ whenever $n \in \Lambda_j$. The system satisfied by $b_j$ is given by:

\begin{equation}\label{toy system v2}
-i\partial_t b_j=-b_j|b_j|^2+2b_{j+1}^2\overline{b_j}+2b_{j-1}^2\overline{b_j}
\end{equation}
with the convention that $b_0(t)=b_{P+1}(t)=0$.

The existence of the frequency set $\Lambda$ is contained in the following proposition:

\begin{proposition}\label{construction of Lambda}
Let $0<s, s\neq 1$. Given parameters $\delta \ll 1$ and $K\gg1$ and a finite set $S$ satisfying the $R-$closure condition, there exists a frequency set $\Lambda=\Lambda_1 \cup \ldots \cup \Lambda_P \subset \Z^2$ for some $P=P(\delta,K)\gg 1$ satisfying Property II$_\Lambda$--Property VIII$_\Lambda$ as well as:

\begin{equation}\label{norm explosion v2}
\frac{\sum_{n \in \Lambda_{P-2}}|n|^{2s}}{\sum_{n\in \Lambda_3}|n|^{2s}} \gtrsim \frac{K^2}{\delta^2}
\end{equation}

In addition, for any $N \geq C(K,\delta, R,S)$ we can make sure that $\Lambda$ consists of $P2^{P-1}$ frequencies $n$ satisfying $N\leq |n|\leq C(\delta, K, S,R)N$.
\end{proposition}

\proof Given $\delta \ll 1, K\gg 1$, we first construct a set $\Lambda^0=\Lambda_1^0 \cup \ldots \cup \Lambda_P^0$ satisfying the $0-$closure condition and Property V$_\Lambda$--Property VIII$_{\Lambda}$, as well as the norm explosion condition $\eqref{norm explosion v2}$. This is exactly the set $\Lambda^0$ constructed\footnote{Strictly speaking, the construction in \cite{CKSTT2} is stated for $s\geq 1$, but an inspection of the proof there (namely pages 102-103) shows that the same holds for $0<s<1$ by relabeling the indices of the generations $\Lambda_j$ as $j \to P-j$. See \cite{CKSTT2} or Proposition 8.2.1 of \cite{Thesis} for details.} in \cite{CKSTT2}.

Choose $v_0 \in \Z^2$ such that for any distinct $n_1,n_2 \in \Lambda^0$ and any distinct $k_1, k_2 \in S$ we have:

$$
v_0.(n_1-n_2) \neq 0 \;\;\;\; \text{ and } v_0.(k_1-k_2)\neq 0
$$

By replacing $R$ by $[R]+1$, we may assume without any loss of generality that $R \in \N$. We will take $\Lambda=N\Lambda^0 -v$ for any $N\geq N_0(K, \delta, S, R)$, where $N_0(K, \delta, S, R)$ and $v$ will be chosen later such that $|v|\ll N_0$. Notice that choosing $N_0>R$ guarantees that if $\Lambda^0$ satisfies the $0-$closure condition then $N\Lambda^0$ satisfies the $R-$closure condition. Denoting
$$
w_4(n_1,n_2,n_3, n)=|n_1|^2-|n_2|^2-|n_3|^2-|n|^2
$$
and noticing that $\omega_4(n_1, \ldots, n)=\omega_4(n_1-w, n_2-w, n_3-w,n-w)$ for any $w \in \Z^2$, we get that $\Lambda=N\Lambda^0-v$ satisfies Property III$_\Lambda$. Property IV$_\Lambda$ is also obviously satisfied if $N_0 > R$ by a similar reasoning. Properties V$_\Lambda$-VIII$_\Lambda$ are trivially inherited by $\Lambda$ from $\Lambda^0$. 

Also, since $N\Lambda^0$ satisfies $\eqref{norm explosion v2}$, then so does $\Lambda$ because $v \ll N\leq |n|$ for any $n \in N\Lambda^0$. Similarly, the last statement in proposition \ref{construction of Lambda} follows easily.

As a result, all what remains is to verify Property II$_\Lambda$. This will specify $v$. We will choose $v=l Rv_0$, where $l \in \N$ and $L\leq l \leq 2L$ for some $C(\delta, K, S) \leq L \ll \frac{N_0}{R|v_0|}$. The exact choice of $l \in [L, 2L]$ will follow from the following counting argument. There are four possible cases to consider:

\underline{\emph{Case 1:}} Suppose that $(n_1,n_2, k)\in \Lambda \times \Lambda \times S$ and $n=n_1-n_2+k$ form a parallelogram. Then $n_1+v, n_2+v, k$ and $n$ also form a parallelogram. Defining $\tilde n_i\in \Lambda^0$ such that $n_i+v=N\tilde n_i$ ($i=1,2$) we have:

\begin{align*}
\omega_4(n_1, n_2,k_3, n)=& |n_1|^2-|n_2|^2+|k|^2-|n|^2\\
=&|n_1+v|^2-|n_2+v|^2+|k|^2-|n|^2-2v(n_1-n_2)\\
=&\omega_4(n_1+v,n_2+v, k, n) -2vN(\tilde n_1-\tilde n_2)
\end{align*}
Note that since $(n_1+v,n_2+v, k)\in N\Lambda^0 \times N\Lambda^0 \times S$, there are $O_{\Lambda^0, S}(1)$ possible values for $\omega_4(n_1+v,n_2+v, k, n)$ (the number of possible tuples in $\Lambda_0\times \Lambda_0\times S$) and these values are independent of $v$. Writing $v=l R v_0$, we directly notice that for each such possible value of $\omega_4(n_1+v,n_2+v, k, n)$, there is at most one choice of $l \in \Z$ that would make $|\omega_4|\leq R$ (since $v_0.(\tilde n_1-\tilde n_2)\neq 0$).

\underline{\emph{Case 2:}} Now suppose that $(n_1, k, n_2) \in \Lambda \times S \times \Lambda$ form a parallelogram with the forth vertex being $n=n_1 -k +n_2$. In this case, $n_1+v, k, n_2+v$ and $n+2v$ also form a parallelogram and:
\begin{align*}
\omega(n_1, k, n_2, n)=&\omega(n_1+v, k+v, n_2+v, n+v)=|n_1+v|^2+|n_2+v|^2-|k+v|^2-|n+v|^2\\
=&|n_1+v|^2+|n_2+v|^2-|k|^2-|n+2v|^2+2v(n+v-k)\\
=&\omega_4(n_1+v, k, n_2+v, k, n+2v)+2v(n+v-k)
\end{align*}
Notice that since $(n_1+v, k , n_2+v)\in N\Lambda_0\times S\times N\Lambda_0$, there are $O_{\Lambda_0, S}(1)$ possible values of $\omega_4(n_1+v, k, n_2+v, k, n+2v)$ and these values are independent of $v$. Writing as before $n_i+v=N\tilde n_i$ with $\tilde n_i\in \Lambda_0$ for $i=1,2$, we have:
\begin{align*}
2v(n+v-k)=&2v\left((n_1+v)+(n_2+v)-v-2k\right)=2Rl v_0\cdot\left(N(\tilde n_1+\tilde n_2)-Rlv_0-2k\right)\\
=&2R \left(Nv_0 \cdot (\tilde n_1+\tilde n_2)l-R|v_0|^2l^2-2k.v_0l\right)\in 2R\Z
\end{align*}
As a result, for each of the $O_{\Lambda_0, K}(1)$ possible values of $\omega_4(n_1+v, k, n_2+v, k, n_2v)$, there are at most two choices of $l$ for which $|\omega(n_1, k, n_2, n)|\leq R$.

\underline{\emph{Case 3:}} Now suppose that $(k_1, k_2, n_3)\in S \times S \times \Lambda$ and $n=k_1-k_2+n_3$ form a parallelogram. We have that:
\begin{align*}
\omega_4(k_1,k_2, n_3, n)=&\omega_4(k_1, k_2, n_3+v, n+v)-2v.(n_3-n)\\
=&\omega_4(k_1, k_2, n_3+v,n+v)-2v(k_1-k_2) \\
=&\omega_4(k_1, k_2, n_3+v,n+v)-2lRv_0.(k_1-k_2).
\end{align*}
As before, there are $O_{\Lambda_0, S}(1)$ possible values for $\omega_4(k_1, k_2, n_3+v,n+v)$ and these values are independent of $v$. For each such value there is at most one choice of $l \in \N$ so that $|\omega_4(k_1,k_2, n_3, n)|\leq R$ (since $v_0\cdot(k_1 -k_2)\neq 0$). 

\underline{\emph{Case 4:}} Finally, suppose that $(k_1, n_2, k_3) \in S\times \Lambda \times S$ and $n=k_1-n_2+k_3$ form a parallelogram. Then,
\begin{align*}
\omega_4(k_1, n_2, k_3, n)&=\omega_4(k_1, n_2+v, k_3, n-v)+2|v|^2+2(n_2+n).v\\
&=\omega_4(k_1, n_2+v, k_3, n-v)+2|v|^2+2(k_1+k_3).v\\
&=\omega_4(k_1, n_2+v, k_3, n-v)+2R\left(R|v_0|^2l^2+2(k_1+k_3).v_0 l\right)
\end{align*}
For each of the $O_{\Lambda_0, S}(1)$ possible values of $\omega_4(k_1, n_2+v, k_3, n-v)$, there are at most two choices of $l$ that can make $|\omega_4(k_1, n_2, k_3, n)|\leq R$. 

In conclusion, for each tuple in $(S\times S \times \Lambda) \cup (\Lambda \times \Lambda \times S)$, there are at most two choices of $l\in \N$ for which a parallelogram determined by this tuple satsifies $|\omega_4|\leq R$. But there are only $C(P, S)=C(\delta, K, S)$ such tuples, so by choosing $L > C(\delta, K, S)$, the pigeonhole principle guarantees the existence of an $L \leq l \leq 2L$ for which all such parallelograms joining $\Lambda$ and $S$ satisfy $|\omega_4|\geq R$. This finishes the proof. \endproof

\begin{remark}\label{remark on construction of Lambda}
It is easy to notice that if $S=\{0\}$ and if $\Lambda^1$ is a set satisfying the conditions of Proposition $\ref{construction of Lambda}$, then so is any dilation $N\Lambda^1$. In particular, one could construct $\Lambda$ in the sublattice $N\Z^2$, which gives the proof of Proposition \ref{construction of Lambda v1}. 
\end{remark}

With Proposition \ref{construction of Lambda} in hand, we are now in a position to finish the proof of Proposition \ref{LT instability c10}. 

Suppose that $\delta \ll 1, K \gg 1$ and that $S$ is a finite set satisfying the $R-$closure condition. Recall that we only need to construct a solution $(a_n(t))$ to $\eqref{NNLS2}$ that is supported in frequency space on a set $\Lambda$ satisfying properties Property I$_\Lambda$--Property VIII$_\Lambda$ and that satisfies $\eqref{a_n(0)}$ and $\eqref{a_n(T)}$.  

Construct $\Lambda$ as in Proposition \ref{construction of Lambda} with $P=P(\delta, K)$ and some $N \geq C(\delta, K, S)$. Using Proposition \ref{diffusion result} we construct a solution to $\eqref{toy system v2}$ satisfying:
\begin{align*}
|b_3(0)|\geq 1-\epsilon&\;\;\;\;|b_j(0)|\leq \epsilon \;\;\;j\neq 3 \\
|b_{P-2}(T)|\geq 1-\epsilon&\;\;\;\; |b_j(T)|\leq \epsilon \;\;\; j\neq P-2
\end{align*}
where $\epsilon=\epsilon(\delta,K)$ is to be specified later and $T=T(\epsilon, P)$.

By construction, defining $\tilde a_n(t)=b_j(t)$ for $n \in \Lambda_j, j=1,\ldots P$ gives a solution to \eqref{RFNLS2}, and hence to $\eqref{NNLS2}$ thanks to Property IV$_\Lambda$. The requested solution $(a_n(t))$ will be a rescaling of $(\tilde a_n)$ that preserves system $\eqref{NNLS2}$, namely:
\begin{equation*}
a_n(t)=\frac{1}{\lambda}\tilde a_n(\frac{t}{\lambda^2}), \quad t\in [0, \lambda^2 T].
\end{equation*}
We choose $\lambda$ so that $\|(a_n(0))\|_{h^s(\Z^2)}\sim \delta$. In fact, since
$$
\|(a_n(0))\|_{h^s(\Z^2)}\sim \frac{1}{\lambda}\left(\sum_{j=1}^P\sum_{n\in \Lambda_j}|n|^{2s}\right)^{1/2}\sim_{\delta, K, S} \lambda^{-1}N^s 
$$
we only need to take $\lambda\sim C(\delta,K,S) N^s$ for $\eqref{a_n(0)}$ to hold. The same calculation as that done in Section \ref{scaling argument} gives that the ratio:

\begin{equation*}
Q=\frac{\|(a_n(\lambda^2T))\|_{h^s(\Z^2)}}{\|(a_n(0))\|_{h^s(\Z^2)}} \gtrsim \frac{K}{\delta}
\end{equation*}
from which $\eqref{a_n(T)}$ follows.
\endproof

\bibliography {bib/network,bib/naming}    
\bibliographystyle {thesis}

\end{document}